\documentclass[12pt,leqno]{article}
\title{\ }
\author{\ }
\begin{document}
\newtheorem{theorem}{Theorem}[section]
\newtheorem{proposition}{Proposition}[section]
\newtheorem{lemma}{Lemma}[section]
\newtheorem{definition}{Definition}[section]
\newtheorem{corollary}{Corollary}[section]
\newtheorem{remark}{Remark}[section]
\newcommand{\bs}{\textcolor{black}}
\newcommand{\iii}{{\, \vert\kern-0.25ex\vert\kern-0.25ex\vert\, }}
\newcommand{\vect}[1]{\mathit{\boldsymbol{#1}}}
\newcommand{\ffi}{\varphi}
\def\Lim{\displaystyle\lim}
\def\Sup{\displaystyle\sup}
\def\Inf{\displaystyle\inf}
\def\Max{\displaystyle\max}
\def\Min{\displaystyle\min}
\def\Sum{\displaystyle\sum}
\def\Frac{\displaystyle\frac}
\def\Int{\displaystyle\int}
\def\n{|\kern -.05cm{|}\kern -.05cm{|}}
\def\Bigcap{\displaystyle\bigcap}
\def\E{{\cal E}}
\def\R{{\bf \hbox{\sc I\hskip -2pt R}}} 
\def\N{{\bf \hbox{\sc I\hskip -2pt N}}} 
\def\Z{{\bf Z}} 
\def\Q{{\bf \hbox{\sc I\hskip -7pt Q}}} 
\def\C{{\bf \hbox{\sc I\hskip -7pt C}}} 
\def\T{{\bf \hbox{\sc I\hskip -5.7pt T}}} 
{\begin{center} {\ } \vskip 2.1cm {\large\bf Stability of Degenerate Schr\"{o}dinger Equation with Harmonic Method
 \\ } \
\\ \ \\

\end{center}
\newcommand{\AL}{\mathcal{A}}
\begin{center}
\begin{tabular}{c}
\sc Khadidja Fekirini, Naima Louhibi \& Abbes Benaissa,
\\ \\
\small  Laboratory of Analysis and Control of PDEs, \qquad \\
\small Djillali Liabes University,\qquad \\
\small P. O. Box 89, Sidi Bel Abbes 22000, ALGERIA.\qquad \\
\small E-mails: fekirini$_{-}$khadidja@yahoo.fr\\
\small     louhibiben@gmail.com\\
\small benaissa$_{-}$abbes@yahoo.com \\
\end{tabular}
\end{center}

\
\\

\begin{abstract}
This paper is devoted to study the well-posedness and stability of degenerate Schr\"{o}dinger equation
with a boundary control acting at the degeneracy. First, we establish the well-posedness of the degenerate
problem $v_t(x,t)+\imath(x^\alpha v_x(x,t))_x=0, \hbox{ with } x \in (0,1)$, controlled by Dirichlet-Neumann
conditions. Then, exponential and polynomial decreasing of the solution are established. This result is optimal and it is
obtained using complex analysis method.
\end{abstract}

\section{Introduction}
This paper is devoted to study the existence and stability of solutions of degenerate Schr\"{o}dinger equation with
a boundary control acting at the degeneracy in suitable sobolev spaces.  More precisely, we consider the following system
\begin{equation}\label{P}
\left\{
\matrix{v_t(x,t)+\imath(x^\alpha v_x(x,t))_x=0\hfill & \hbox{ in } (0,1)\times(0,+\infty),\hfill & \cr
(x^\alpha v_x)(0,t)=\imath\rho \partial^{\tilde{\alpha}, \eta} v(0,t)\hfill &  \hbox{ on } (0,+\infty),\hfill & \cr
v_x(1,t)=0\hfill & \hbox{ on } (0,+\infty),\hfill & \cr
v(x,0)=v_0(x)\hfill & \hbox{ on } (0,1),\hfill & \cr}\right.
\end{equation}
where $0<\alpha<1, \rho>0$ and the term
$\partial^{\tilde{\alpha}, \eta}$ stands for the generalized fractional integral of order $0<\tilde{\alpha}\leq 1$
(see {\bf\cite{choi}}), which is given by
$$
\partial^{\tilde{\alpha}, \eta}w({t})=\left\{\matrix{w(t)\hfill &\hbox{ for } \tilde{\alpha}=1,\ \eta\geq 0,\hfill \cr
\Frac{1}{\Gamma(1-\tilde{\alpha})}\Int_{0}^{{t}} ({t}-s)^{-\tilde{\alpha}}e^{-{\eta}({t}-s)}w(s)\, ds,\hfill & \hbox{ for } 0<\tilde{\alpha}<1,\ \eta\geq 0.\hfill \cr}\right.
\label{cap}
$$
The controllability and stabilization of Schr\"{o}dinger equations without degeneracies have attracted
considerable attention over the past years. Under the so-called geometric control condition, it is proved
by G. Lebeau {\bf\cite{lebeau}} that the Schr\"{o}dinger equation is exactly controllable for arbitrary short time.

In {\bf\cite{M}}, Machtyngier addressed the exact
controllability in $H^{-1}(\Omega)$, with $\Omega$ is a bounded domain, where Dirichlet boundary condition
in $L^{2}(\Omega)$. The approaches adapted are HUM (Hilbert Uniqueness Method) and multipliers techniques.

The boundary stabilization of the Schr\"{o}dinger equation has also received a lot of attention. For an introduction,
see {\bf\cite{Ire-Trig}}, where Lasiecka and Triggiani examine solution existence, uniqueness, and uniform boundary
stability at the energy level in $ L^{2}(\Omega)$ for the n-dimensional linear Schr\"{o}dinger equation within
a bounded open domain. This system is given by:
$$
\left\{ \begin{array}{lr}
        u_{t}+i\Delta u=0&\Omega\times(0,+\infty),\\
        \frac{\partial u }{\partial \nu  }=iu &x\in \Gamma_1, t \geq 0, \\
        u= 0 \ &x\in \Gamma_2   , t \geq0, \\
        u(x,0)=u_{0}(x),u_{t}(x,0)=u_{1}(x)&x\in\Omega.
    \end{array} \right.
    $$
The authors adopted semigroup theory to show the global existence of the system and thereafter determined
an optimal decay result using the multiplier method.

A similar study was accomplished in {\bf\cite{My-Zau}},
where $iu$ was replaced with $im(x)u_t$. In this case, the authors proved exponential decay in both the
$L^{2}$-norm and the $H^{1}$-norm by employing the same approach while imposing geometric control conditions on
the boundary.

In {\bf\cite{KGS}}, the problem treated is the following
$$
\quad \left \{
\begin{array}{ll}
\Phi_t(x,t)=-\imath\Phi_{xx}(x,t), & 0<x<1,\ \  t>0, \\
\Phi(1,t)= u(t),  & t\geq 0,\\
\Phi_x(0,t)=0, &  t\geq 0,\\
y(t)=\Phi_x(1,t), &  t\geq 0,
\end{array}
\right.
$$
where $u(t)$ is the "input" and $y(t)$ is the "output".
Two novel control designs are proposed to exponentially stabilize the system.
Furthermore, S. Nicaise and S. Rebiai {\bf\cite{Rebiai}} examined the influence of time delays on the boundary and
internal feedback stabilization of the multidimensional Schr\"{o}dinger equation, which is usually used to model the
behavior of quantum systems. Their study centred on how these time delays influence the stabilization process, focusing to
provide insights into the dynamics and control of such systems across different contexts. The system is given by:
$$
\left\{\begin{array}{lr}
y_{t}(x,t)+i\Delta y(x,t)  = 0,& x \in  \Omega, t> 0,\\
y(x,0)  =y_{0}(x,t)&x \in  \Omega, t> 0,\\
y(x,t)= 0&x \in  \Gamma_{0}, t> 0,\\
\frac{\partial y}{\partial \nu}=i\mu_{1}y(x,t)+i\mu_{2} y(x,t-\tau),&\in  \Gamma_{1}, t> 0,\\
y(x,t-\tau)=f_{0}(x,t-\tau)&x \in \Omega,0<t<\tau.\\
\end{array}\right.
$$
On the contrary, when the principal part is degenerate
not much is known in the literature, despite that many problems that are pertinent for applications
are modeled by Schr\"{o}dinger equations degenerating at the boundary
of the space domain.

\noindent
In {\bf\cite{jawad}}, the authors considered the following Schr\"{o}dinger equation
\begin{equation}\label{zaz}
\left\{
\matrix{v_t(x,t)+\imath(x^\alpha v_x(x,t))_x=0\hfill & \hbox{ in } (0,1)\times(0,+\infty),\hfill & \cr
v(0,t)=0\hfill &  \hbox{ on } (0,+\infty),\hfill & \cr
v_t(1,t)+v_x(1,t)+v(1,t)=0\hfill & \hbox{ on } (0,+\infty),\hfill & \cr
v(x,0)=v_0(x)\hfill & \hbox{ on } (0,1),\hfill & \cr}\right.
\end{equation}
They proved that the solution decays exponentially in an appropriate energy space.
Moreover, the degeneracy does not affect the decay rates of the energy.

Here, the situation is different since we impose a damping at point $x=0$, where the degeneracy of the elliptic
operator $(x^{\alpha}\partial_x v)_x$ holds, which turns out to be a more challenging issue.

To our best knowledge, this is the first attempt to study the global decaying solutions for a degenerate Schr\"{o}dinger
equation under a control acting on the degenerate boundary.
Moreover, the enegy method based on multiplier techniques used in {\bf\cite{jawad}} do not seem to be work in the case of
a feedback acting at a degenerate point $x=0$.

In this work, we are interested in studying precisely this issue, extending the results obtained in {\bf\cite{ZBA}},
where the authors discuss the same issue in the case of wave equations.
We obtain new results on decay estimates depending on parameters $\alpha$ and $\tilde{\alpha}$.

This paper is organised as follows. In section 2, we give some preliminaries. In section 3,
the well-posedness results of the system $(\ref{P})$ are given using semigroup theory. In section 4,
we prove an asymptotic and polynomial decay using Borichev-Tomilov Theorem. In section 5, we prove lack of
exponential stability using Rouch\'{e}'s Theorem.
of the obtained system $(\ref{P})$ and we prove an optimal decay rate.

\section{Preliminaries}
In this section, we introduce notations, definitions and
propositions that will be used later. First we introduce some  weighted Sobolev spaces:
$$
H^1_{\alpha}(0,1)=\left\{v \in L^2(0,1), v \hbox{  is locally absolutely continuous in }(0,1], \ x^{\alpha /2}v_x\in L^2(0,1)\right\},
$$
and
$$
H^2_{\alpha}(0,1)=\left\{v \in L^2(0,1), v\in H^1_{\alpha}(0,1), \ x^{\alpha}v_x\in H^1(0,1)\right\},
$$
where $H^1(0,1)$ represent the classical Sobolev space.\\
We remark that $H^1_{\alpha}(0,1)$ is a Hilbert space with the scalar product
$$
(u, v)_{H_\alpha^1(0, 1)}=\Int_{0}^{1}(u\overline{v}+ x^{\alpha} u'(x)\overline{v'(x)})\, dx,\quad\forall u, v\in H_\alpha^1(0, 1).
$$
\begin{remark}
we have that
$$
\matrix{v^2(1)&=&\Int_{0}^{1}(x v^2(x))_x\, dx=\Int_{0}^{1}(v^2(x)+2 xv(x)v'(x))\, dx\hfill \cr
&\leq & 2\Int_{0}^{1}v^2(x)\, dx+\Int_{0}^{1}x^2 (v'(x))^2\, dx\hfill \cr
&\leq & 2\Int_{0}^{1}v^2(x)\, dx+\Int_{0}^{1}x^{\alpha} (v'(x))^2\, dx\leq 2 \|v\|_{H_{\alpha}^1(0, 1)}^2.\hfill \cr}
$$
Moreover, we have
$$
\matrix{|v(x)|&=&\left|-\Int_{x}^{1}v'(x)\, dx+v(1)\right|\hfill \cr
&\leq & \Int_{0}^{1}|v'(x)|\, dx+|v(1)|\hfill \cr
&\leq & \left(\Int_{0}^{1}\frac{1}{x^{\alpha}}\, dx\right)^{1/2}\left(\Int_{0}^{1}x^{\alpha}|v'(x)|^2\, dx\right)^{1/2}+|v(1)|\hfill \cr
&\leq & \left(\frac{1}{\sqrt{1-\alpha}}+\sqrt{2}\right)\|v\|_{H_{\alpha}^1(0, 1)}.\hfill \cr}
$$
Hence
$$
H_{\alpha}^1(0, 1)\hookrightarrow C([0, 1]).
$$

\end{remark}

\subsection{Augmented model}
In this section we reformulate $(\ref{P})$ into an augmented system. For that, we need the following proposition.
\begin{proposition}[see {\bf\cite{mbod}}]
Let $\mu$ be the function:
\begin{equation}
\mu(\xi)=|\xi|^{(2\tilde{\alpha}-1)/2},\quad -\infty<\xi<+\infty,\
0<\tilde{\alpha}<1. \label{e1}
\end{equation}
Then the relationship between the `input' U and the `output' O of
the system
\begin{equation}
\partial_t\phi(\xi, t)+(\xi^{2}+\eta)\phi(\xi, t) -U(t)\mu(\xi)=0,\quad -\infty<\xi<+\infty,\eta\geq 0, t> 0,
\label{e6}
\end{equation}
\begin{equation}
\phi(\xi, 0)=0, \label{e7}
\end{equation}
\begin{equation}
O(t)=(\pi)^{-1}\sin(\tilde{\alpha}\pi)\Int_{-\infty}^{+\infty}\mu(\xi)\phi(\xi, t)\, d\xi,
\label{e8}
\end{equation}
where $U\in C^0([0,+\infty))$, is given by
\begin{equation}
O=I^{1-\tilde{\alpha}, \eta}U, \label{e9}
\end{equation}
where
$$
[I^{\alpha,
\eta}f](t)=\Frac{1}{\Gamma(\tilde{\alpha})}\Int_{0}^{t}(t-\tau)^{\tilde{\alpha}-1}e^{-\eta(t-\tau)}f(\tau)\, d\tau.
$$
\label{th2}
\end{proposition}
\begin{lemma}[see {\bf\cite{achour.1}}]
If $\lambda\in D_{\eta}=\C\backslash]-\infty, -\eta]$ then
$$
\Int_{-\infty}^{+\infty}\Frac{\mu^2(\xi)}{\lambda+\eta+\xi^2}\, d\xi
=\Frac{\pi}{\sin\tilde{\alpha}\pi}(\lambda+\eta)^{\tilde{\alpha}-1}.
$$
\label{achour}
\end{lemma}
Using now Proposition \ref{th2} and relation (\ref{e9}), system $(\ref{P})$ may be recast into the following
augmented system
$$
\left\{\matrix{v_t(x,t)+\imath(x^\alpha v_x(x,t))_x=0,\hfill &\cr
\phi_t(\xi, t)+(\xi^{2}+\eta)\phi(\xi, t) -v(0, t)\mu(\xi)=0,\hfill & -\infty <\xi<+\infty,\ \ \hfill t>0,\hfill&\cr
(x^\alpha v_x)(0,t)=i\zeta \Int_{-\infty}^{+\infty}\mu(\xi)\phi(\xi, t)\, d\xi,
\hfill & \hfill \cr
v_x(1,t)=0, \hfill & \hfill \cr
v(x,0)=v_0(x), \phi(\xi, 0)=0,\hfill &\cr }\right.
\leqno{(P')}
$$
where $\zeta=\varrho(\pi)^{-1}\sin(\tilde{\alpha}\pi)$.

We define the energy associated with the solution of the problem $(\ref{P})$ by
\begin{equation}\label{27}
E(t)=\frac{1}{2}\int_{0}^{1}|v(x,t)|^2dx+\frac{\zeta}{2}\Int_{-\infty}^{+\infty}|\phi(\xi, t)|^2\, d\xi.
\end{equation}

\begin{lemma}
The energy functional defined by $(\ref{27})$ decays as follows
$$
E'(t)=-\zeta \Int_{-\infty}^{+\infty}(\xi^{2}+\eta)|\phi(\xi,t)|^2d\xi\leq 0.
$$
\end{lemma}
{\bf Proof.}
Multiplying the equation of $(\ref{10})$ by $\overline{v}$, integrating over $(0,1)$, applying integration by parts and
using boundary conditions we obtain
\begin{equation}\label{ep338}
\int_{0}^{1}v_t\overline{v} \ dx =- \zeta \Int_{\R}\mu(\xi)  \phi(\xi)d\xi\overline{v}(0, t)
+\imath \int_{0}^{1}x^\alpha |v_x|^2 \ dx,
\end{equation}
Multiplying the second equation in (\ref{P}) by $\zeta \overline{\phi}$ and integrating over $(-\infty,+\infty)$, we get:
\begin{equation}\label{ep3}
\ \ \ \ \Frac{\zeta}{2}\Frac{d}{dt} \|\phi\|^{2}_{L^2(\R)}
+ \zeta \Int_{\R}(\xi^{2}+\eta)||\phi(\xi,t)|^2d\xi
- \zeta \Re\Int_{\R}\mu(\xi)  \ \overline{\phi}(\xi,t)d\xi v(0, t) = 0.
\end{equation}
We take the sum of (\ref{ep3}) and the real part of (\ref{ep338}), we get
$$
E'(t)=-\zeta \Int_{\R}(\xi^{2}+\eta)||\phi(\xi,t)|^2d\xi\leq 0.
$$
\begin{remark}
In the case $\tilde{\alpha}=1$, we take $\varrho v(0, t)$ instead of $\varrho\partial^{\alpha, \eta} v(0, t)$.
We do not need to introduce an augmented system. In this case the operator ${\cal A}$ takes the form
\begin{equation}
{\tilde{\cal A}}v= -i (x^\alpha u_x)_x
 \label{eqOP}
\end{equation}
with domain
\begin{equation}
D({\tilde{\cal A}})=\left\{\matrix{
u\in H_{\alpha}^{2}(0,1),\hfill
\cr ( x^\alpha u_x)(0)=i\varrho v(0)=0, v_x(1)=0 \hfill \cr }\right\},
\label{eqOP1}
\end{equation}
where
$$
{\tilde{\cal H}}=L^{2}(0,1)
$$
with inner product
$$
\left\langle  v,\tilde{v}\right\rangle_{{\tilde{\cal H}}}=\Int_{0}^{1}v \ {\overline {\tilde v}}dx.
$$
The well-posedness result follows exactly as in the case $0< \tilde{\alpha} <1$. Moreover, the energy function is defined as
\begin{equation}
{\tilde E}(t)=\Frac{1}{2}\int_{0}^{1}|v|^2\, dx
\label{eqOP2}
\end{equation}
and decays as follows
$$
{\tilde E}'(t)=-\varrho|v(0, t)|^2\leq 0.
$$
\end{remark}

\hfill$\diamondsuit$\\

\section{Well-posedness}

\quad This section is concerned to the well-posedness results of the problem $(\ref{P})$ using a semigroup approach and the
Lumer-Philips Theorem.\\

We introduce the Hilbert space
$$
{\mathcal{H}}=L^{2}(0,1)\times L^2(\R),
$$
with the following inner product
$$
\langle V,\widetilde{V} \rangle_{\mathcal{H}}=\int_{0}^{1}v(x)\overline{\widetilde{v}}(x)dx+\zeta\Int_{-\infty}^{+\infty}\phi(\xi)\overline{\widetilde{\phi}}(\xi)\, d\xi
$$
for all $V, \widetilde{V} \in {\mathcal{H}}$ with $V=(v, \phi)^{T}$ and $\widetilde{V}=(\widetilde{v}, \widetilde{\phi})^{T}$.
The problem $(\ref{P})$ can be written as
\begin{equation}\label{10}
 \quad \left \{
\begin{array}{ll}
V_t={\AL}V,\\
V(0)=V_0,
\end{array}
\right.
\end{equation}
where the operator ${\mathcal{A}}$
is defined by
$$
{\cal A}V=\pmatrix{-\imath (x^\alpha v_x)_x  \cr
-(\xi^2+\eta)\phi+\mu(\xi)v(0) \cr}
$$
with domain
$$
D({\mathcal{A}})=\left\{\matrix{ (v, \phi)\in {\cal H}: v\in H_{\alpha}^{2}(0, 1),  v_x(1)=0, \ (x^\alpha v_x)(0)=\imath\zeta\Int_{-\infty}^{\infty}\mu(\xi)\phi(\xi)\, d\xi & \cr
-(\xi^2+\eta)\phi+\mu(\xi)v(0)\in L^2(\R), |\xi|\phi\in L^2(\R)& \cr}\right\}.
$$
We will show that the operator ${\mathcal{A}}$ generates a $C_0$-semigroup of contractions in the Hilbert space
${\mathcal{H}}$.

First, we prove that the operator ${\mathcal{A}}$ is dissipative. We have for every $V\in D({\mathcal{A}})$,
\begin{equation}
\Re\langle{\mathcal{A}}V, V\rangle_{{\mathcal{H}}}=
-\zeta \int_{-\infty}^{+\infty}(\xi^{2}+\eta)|\phi(\xi)|^{2}\, d\xi\leq 0.
\label{zz18}
\end{equation}
Next, we prove that the operator $(\lambda I-{\mathcal{A}})$ is surjective for $\lambda>0$. For this purpose, let us
take $F\in {\mathcal{H}}$, we search $V\in D({\mathcal{A}})$ such that
\begin{equation}\label{14}
(\lambda I-{\mathcal{A}})V=F.
\end{equation}
From equation (\ref{14}), we get the following system of equations
\begin{equation}
\left\{\matrix{\imath\lambda v-(x^\alpha v_x)_x=\imath f_1,  \hfill & \cr
\lambda \phi+(\xi^{2}+\eta)\phi-v(0)\mu(\xi)=f_2.\hfill & \cr}\right.
\label{e18}
\end{equation}
By $(\ref{e18})_2$ we can find $\phi$ as
\begin{equation}
\phi(\xi)=\Frac{v(0)\mu(\xi)+f_2(\xi)}{\lambda +\xi^{2}+\eta}.
\label{m18}
\end{equation}
Solving $(\ref{e18})_1$ is equivalent to finding $v\in H^2_{\alpha}(0,1)$
such that, for all $w \in H^1_{\alpha}(0,1)$
\begin{equation}\label{19}
\int_{0}^{1}\imath\lambda v \overline{w}\ dx-\int_{0}^{1}(x^\alpha v_x)_x\overline{w}\ dx
=\int_{0}^{1}\imath f_1\overline{w}\ dx.
\end{equation}
Using $(\ref{19})$, the boundary conditions and (\ref{m18}), the function $v$ satisfies the following equation, for all
$w \in H^1_{\alpha}(0,1)$
$$
\imath\int_{0}^{1}\lambda v \overline{w}\ dx+\imath\rho (\lambda+\eta)^{\tilde{\alpha}-1}v(0)\overline{w}(0)+\int_{0}^{1}x^\alpha v_x\overline{w}_x\ dx
=\imath\int_{0}^{1} f_1\overline{w}\ dx-i\zeta\Int_{-\infty}^{+\infty}\Frac{\mu(\xi)f_2(\xi)}{\lambda +\xi^{2}+\eta}\, d\xi \overline{w}(0),
$$
Multiplying this equation by $(1-i)$, we obtain
$$
\matrix{(1-i)\int_{0}^{1}(i\lambda v \overline{w}+\int_{0}^{1}x^\alpha v_x\overline{w}_x)\ dx
+(1-i)\imath\rho (\lambda+\eta)^{\tilde{\alpha}-1}v(0)\overline{w}(0)= \hfill &\cr
(1-i)\imath\int_{0}^{1} f_1\overline{w}\ dx
-i(1-i)\zeta\Int_{-\infty}^{+\infty}\Frac{\mu(\xi)f_2(\xi)}{\lambda +\xi^{2}+\eta}\, d\xi \overline{w}(0), &\cr}
$$
which is of the form
\begin{equation}\label{24}
{\mathbf{b}}(v,w)={\mathbf{l}}(w),
\end{equation}
where ${\mathbf{b}}:[H^1_{\alpha}(0,1)\times H^1_{\alpha}(0,1)]\longrightarrow \C$
is the sesquilinear form defined by
$$
{\mathbf{b}}(v,w)=(1-i)\imath\int_{0}^{1}\lambda v \overline{w}\ dx+(1-i)\int_{0}^{1}x^\alpha v_x\overline{w}_x\ dx
+(1-i)\imath\rho (\lambda+\eta)^{\tilde{\alpha}-1}v(0)\overline{w}(0),
$$
and ${\mathbf{l}}:H^1_{\alpha}(0,1)\longrightarrow \C$
is the antilinear form given by
$$
{\mathbf{l}}(w)=(1-i)\imath\int_{0}^{1} f_1\overline{w}\ dx
-i(1-i)\zeta\Int_{-\infty}^{+\infty}\Frac{\mu(\xi)f_2(\xi)}{\lambda +\xi^{2}+\eta}\, d\xi \overline{w}(0).
$$
It is easy to verify that ${\mathbf{b}}$ is continuous and coercive and ${\mathbf{l}}$ is continuous, therefore
using the Lax-Milgram theorem, we conclude that the problem $(\ref{24})$ admits a
unique solution $v\in H^1_{\alpha}(0,1)$, for all $\lambda>0$.
Now, if we consider $w\in {\cal D}(0, 1)$  in (\ref{24}), then $v$ solves in ${\cal D}'(0, 1)$
$$
\lambda v+i(x^\alpha v_x)_x= f_1
$$
and thus $(x^\alpha v_x)_x\in L^2(\Omega)$.

Using Green's formula in (\ref{24}), we get
$$
[(x^\alpha v_x){\overline w}]_0^1=
-i\zeta\Int_{-\infty}^{+\infty}\Frac{\mu(\xi)f_2(\xi)}{\lambda +\xi^{2}+\eta}\, d\xi \overline{w}(0)
-\imath\rho (\lambda+\eta)^{\tilde{\alpha}-1}v(0)\overline{w}(0).
$$
Consequently, defining $\phi$ by $(\ref{m18})$, we conclude that
\begin{equation}
[(x^\alpha v_x){\overline w}]_0^1=-i\zeta\Int_{-\infty}^{+\infty}\mu(\xi)\phi(\xi)\, d\xi \overline{w}(0).
\label{gghh}
\end{equation}
If we take $w(x)=x$, we find $v_x(1)=0$. If we take $w(x)=1$, we find
$$
(x^\alpha v_x)(0)=i\zeta\Int_{-\infty}^{+\infty}\mu(\xi)\phi(\xi)\, d\xi.
$$
In order to achieve the existence of $V\in D({\cal A})$,  we require to prove $\phi, |\xi|\phi\in L^{2}(\R)$
$$
\Int_{\R}|\phi(\xi)|^2\, d\xi\leq 2\Int_{\R}\Frac{|f_2(\xi)|^2}{(\xi^{2}+\eta+\lambda)^2}  \, d\xi
+2 |v(0)|^2\Int_{\R}\Frac{\mu(\xi)^2}{(\xi^{2}+\eta+\lambda)^2}  \, d\xi.
$$
On the other hand, using the fact that $f_2\in L^2(\R)$, we get
$$
\Int_{\R}\Frac{|f_2(x, \xi)|^2}{(\xi^{2}+\eta+\lambda)^2}  \, d\xi\leq
\Frac{1}{(\eta+\lambda)^2}\Int_{\R}|f_2(\xi)|^2  \, d\xi< +\infty.
$$
and
$$
\Int_{\R}\Frac{\mu(\xi)^2}{(\xi^{2}+\eta+\lambda)^2}  \, d\xi\leq \Frac{1}{(\eta+\lambda)}\Int_{\R}\Frac{\mu(\xi)^2}{\xi^{2}+\eta+\lambda}  \, d\xi< +\infty.
$$
Thus $\phi \in L^{2}(\R)$. Next, using again (\ref{m18}), we get
$$
\Int_{\R}|\xi\phi(\xi)|^2\, d\xi\leq 2\Int_{\R}\Frac{|\xi|^2|f_2(\xi)|^2}{(\xi^{2}+\eta+\lambda)^2}  \, d\xi
+2 |v(0)|^2\Int_{\R}\Frac{|\xi|^2\mu(\xi)^2}{(\xi^{2}+\eta+\lambda)^2}  \, d\xi.
$$
Using the fact that $f_2\in L^2(\R)$, we obtain
$$
\Int_{\R}\Frac{|\xi|^2|f_2(\xi)|^2}{(\xi^{2}+\eta+\lambda)^2}  \, d\xi\leq
\Frac{1}{(\eta+\lambda)}\Int_{\R}|f_2(\xi)|^2  \, d\xi< +\infty.
$$
and
$$
\Int_{\R}\Frac{|\xi|^2\mu(\xi)^2}{(\xi^{2}+\eta+\lambda)^2}  \, d\xi\leq
\Int_{\R}\Frac{\mu(\xi)^2}{\xi^{2}+\eta+\lambda}  \, d\xi< +\infty.
$$
Thus $\xi\phi \in L^{2}(\R)$. Moreover
$$
-(\xi^{2} + \eta)\phi+\mu(\xi)v(0)=\lambda\phi-f_{2}\in L^{2}(\R).
$$
So applying the Hille-Yoshida Theorem we have the following result.
\begin{theorem}(Existence and uniqueness)
\begin{enumerate}
\item If $V_0 \in D({\mathcal{A}})$, then the problem $(\ref{10})$ has a unique strong solution
$$
V\in {\mathcal{C}}^0({\R}_+,D({\mathcal{A}}))\cap {\mathcal{C}}^1(\R_+,\mathcal{H}).
$$
\item If $V_0 \in \mathcal{H}$, then the problem $(\ref{10})$ has a unique weak solution
$$
V\in {\mathcal{C}}^0({\R}_+,\mathcal{H}).
$$
\end{enumerate}
\end{theorem}

\hfill$\diamondsuit$\\

\section{Asymptotic and decay estimates of solutions}
In this section, we will study the stability of solution associated with the problem $(\ref{10})$,
for this purpose we need the following theorem.
\begin{lemma}\label{5.1}{\bf\cite{AB}}
Let ${\mathcal{A}}$ be the generator of a uniformly bounded $C_0$-semigroup $\{S(t)\}_{t\geq0}$ on a Hilbert space
${\mathcal{H}}$. If
\begin{enumerate}
  \item ${\mathcal{A}}$ does not have eigenvalues on $\imath{\R}$.
  \item The intersection of the spectrum $\sigma({\mathcal{A}})$ with $\imath{\R}$ is at most a countable set,
  \end{enumerate}
  then the semigroup ${S(t)}_{t\geq0}$ is asymptotically stable, i.e. $\|S(t)z\|_{\mathcal{H}}\rightarrow0$ when $t\rightarrow\infty$
  for any $z\in{\mathcal{H}}$.
\end{lemma}

We will use this theorem to prove the strong stability of the $C_0$-semigroup $e^{t{\mathcal{A}}}$.
Our main result is the following theorem.
\begin{theorem}\label{5.3}
The $C_0$-semigroup $e^{t{\mathcal{A}}}$ is strongly stable in ${\mathcal{H}}$, i.e. for all $V_0 \in {\mathcal{H}}$,
the solution of the problem $(\ref{10})$ verify
$$
\lim_{t\rightarrow\infty}\|e^{t{\mathcal{A}}}V_0\|_{\mathcal{H}}=0.
$$
\end{theorem}

In order to prove the Theorem $\ref{5.3}$, we need the following  two lemmas.
\begin{lemma}\label{l5.1}
The operator ${\mathcal{A}}$ does not have eigenvalues on $\imath{\R}$.
\end{lemma}
{\bf Proof.}
We start with the first case $\imath \lambda =0$.
The equation ${\mathcal{A}}V=0$ leads to
$$
\left\{\matrix{-(x^\alpha v_x)_x=\hfill &\cr
-(\xi^2+\eta)\phi+v(0)\mu(\xi)=0\hfill &\cr}\right.
$$
From (\ref{zz18}), we have $\phi\equiv 0$ and then
\begin{equation}
\quad \left\{
\begin{array}{ll}
v(0)=0,\\
(x^\alpha v_x)(0)=0\\
\end{array}
\right.
\label{zz22}
\end{equation}
with $v_x(1)=0$. Hence
\begin{equation}\label{*}
(x^\alpha v_x)(x)=0,
\end{equation}
So, for all $ x \in (0,1)$ $$v_x(x)=0,$$
then $v(x)=\widetilde{c}$, where $\widetilde{c}$ is a constant, as $v(0)=0$, so $$v=0.$$
Hence $\imath\lambda=0$ is not an eigenvalue of ${\mathcal{A}}$.\\

Next, we study the case $\imath\lambda\neq0$. Let us suppose that $\lambda \in {\R}-\{0\}$ such that
$${\mathcal{A}}V=\imath\lambda V,$$ with $V\neq 0$, then we get
$$
\left\{\matrix{(x^\alpha v_x)_x=-\lambda v \hfill &\cr
-(\xi^2+\eta)\phi+v(0)\mu(\xi)=i\lambda\phi\hfill &\cr}\right.
$$
Using (\ref{zz18}), we get $\phi\equiv 0$.
so, we obtain the following system
\begin{equation}\label{54}
\quad \left \{
\begin{array}{ll}
\lambda v+(x^\alpha v_x)_x=0, \quad \hbox{ on }(0,1),\\
(x^\alpha v_x)(0)=v(0)=0,\\
v_x(1)=0.
\end{array}
\right.
\end{equation}
This type of problems can be solved using the Bessel functions.
The solution of the $(\ref{54})_1$ is given by
$$
v(x)=c_1\theta_+(x)+c_2\theta_-(x),
$$
where $c_1$ and $c_2$ are two constants, and $\theta_+$ and $\theta_-$ are defined by
\begin{equation}\label{theta}
\theta_+(x)=x^{\frac{1-\alpha}{2}}J_{\nu_{\alpha}}\left(\frac{2}{2-\alpha}\mu x^{\frac{2-\alpha}{2}}\right)
\quad \hbox{ and  } \
\theta_-(x)=x^{\frac{1-\alpha}{2}}J_{-\nu_{\alpha}}\left(\frac{2}{2-\alpha}\mu x^{\frac{2-\alpha}{2}}\right)
\end{equation}
with $\mu=\sqrt{\lambda}$, $\nu_\alpha=\frac{1-\alpha}{2-\alpha}$,\\
$$
J_{\nu_{\alpha}}(y)=\sum_{m=0}^{\infty}\frac{(-1)^m}{m!\Gamma (m+\nu_{\alpha}+1)}\left(\frac{y}{2}\right)^{2m+\nu_{\alpha}}=\sum_{m=0}^{\infty}c^+_{\nu,m}y^{2m+\nu_{\alpha}}
$$
and
$$
J_{-\nu_{\alpha}}(y)=\sum_{m=0}^{\infty}\frac{(-1)^m}{m!\Gamma (m-\nu_{\alpha}+1)}\left(\frac{y}{2}\right)^{2m-\nu_{\alpha}}=\sum_{m=0}^{\infty}c^-_{\nu,m}y^{2m-\nu_{\alpha}},
$$
$J_{\nu_{\alpha}}$ and $J_{-\nu_{\alpha}}$ are Bessel functions of the first kind of order $\nu_{\alpha}$ and $-\nu_{\alpha}$.\\

We can verify that $\theta_+$ and $\theta_-$ $\in H^1_{\alpha}(0,1)$, indeed, in the neighborhood of zero we have
$$
\theta_+(x)\sim d^+x^{1-\alpha}, \qquad
x^{\alpha/2}\ \theta'_+(x)\sim(1-\alpha) d^+x^{-\alpha/2},
$$
$$
\theta_-(x)\sim d^-,\qquad \qquad
x^{\alpha/2}\ \theta'_-(x)\sim(2-\alpha) d^-x^{1-\alpha/2},
$$
where
\begin{equation}\label{dd}
d^+=c^+_{\nu_{\alpha,0}}\left(\frac{2}{2-\alpha}\mu\right)^{\nu_{\alpha}} \ \hbox{ and \ \ }
d^-=c^-_{\nu_{\alpha,0}}\left(\frac{2}{2-\alpha}\mu\right)^{{-\nu_{\alpha}}}.
\end{equation}

Then the condition $(\ref{54})_2$ become
$$
\left\{\matrix{c_1(1-\alpha)d^+&=&0, \hfill \cr
c_2d^-=0,\hfill \cr}\right.
$$
while $(\ref{54})_3$ become
$$
c_1\theta_+'(1)+c_2\theta_-'(1)=0.
$$
Hence $v=0$.
Therefore $V = 0$, which contradicts $\|V\|_{\cal H} = 1$. This completes the proof of Lemma \ref{l5.1}.

\hfill$\diamondsuit$\\

\begin{lemma}
(a) If $\eta=0$, then the operator $(\imath\lambda I-{\mathcal{A}})$ is surjective for any real number $\lambda\not=0$.\\
(b) If $\eta>0$, then $(\imath\lambda I-{\mathcal{A}})$ is surjective for any $\lambda\in \R$.
\label{l5.2}
\end{lemma}
{\bf Proof.}
We will examine two cases.\\\\
\textbf{Case $1$}: $\lambda\neq0$.\\
Let $F\in {\mathcal{H}}$ be given and let $V\in D({\mathcal{A}})$ be such that
\begin{equation}\label{56}
(\imath\lambda-{\mathcal{A}})V=F,
\end{equation}
so, we have
\begin{equation}\label{58}
\left\{\matrix{-\lambda v-(x^\alpha v_x)_x=\imath f_1,  \hfill & \cr
i\lambda \phi+(\xi^{2}+\eta)\phi-v(0)\mu(\xi)=f_2,\hfill & \cr}\right.
\end{equation}
together with the conditions
\begin{equation}\label{cond}
\quad \left \{
\begin{array}{ll}
(x^\alpha v_x)(0)=i\zeta\Int_{-\infty}^{+\infty}\mu(\xi)\phi(\xi)\, d\xi,\\
v_x(1)=0.
\end{array}
\right.
\end{equation}
From $(\ref{58})_2$ and (\ref{cond}), we get
\begin{equation}
(x^{\alpha}v_x)(0)=i\varrho (i\lambda+\eta)^{\tilde{\alpha}-1} v(0)
+i\zeta\Int_{-\infty}^{+\infty}\Frac{\mu(\xi)f_2(\xi)}{i\lambda+\xi^{2}+\eta}\, d\xi.
\label{cond22}
\end{equation}
Suppose that $v$ is a solution of $(\ref{58})_1$, so the function $\Psi$ defined by
\begin{equation}\label{70z}
v(x)=x^{\frac{1-\alpha}{2}}\Psi\left(\frac{2}{2-\alpha} \mu
x^{\frac{2-\alpha}{2}} \right),
\end{equation}
with $\mu=\sqrt{\lambda}$, is a solution of the following inhomogeneous Bessel equation
$$
y^2\Psi''(y)+y\Psi'(y)+\left(y^2-\left( \frac{\alpha-1}{2-\alpha}\right) ^2 \right) \Psi(y)=
-\left( \frac{2}{2-\alpha}\right)^2\left( \frac{2-\alpha}{2\mu}y\right) ^{\frac{3-\alpha}{2-\alpha}}
\imath f_1\left(\left( \frac{2-\alpha}{2\mu}y\right) ^{\frac{2}{2-\alpha}} \right).
$$
We can write $\Psi$ as
\begin{equation}\label{70'z}
\Psi(y)=AJ_{\nu_\alpha}(y)+BJ_{-\nu_\alpha}(y)-\frac{\pi}{2\sin \nu_\alpha \pi}\int_{0}^{y}\frac{\widetilde{f}(s)}{s}(J_{\nu_\alpha}(s)J_{-\nu_\alpha}(y)-J_{\nu_\alpha}(y)J_{-\nu_\alpha}(s))ds,
\end{equation}
where
$$
\widetilde{f}(s)=-\left( \frac{2}{2-\alpha}\right)^2\left( \frac{2-\alpha}{2\mu}s\right) ^{\frac{3-\alpha}{2-\alpha}}
\imath f_1\left(\left( \frac{2-\alpha}{2\mu}s\right) ^{\frac{2}{2-\alpha}} \right).
$$
Using $(\ref{theta})$, $(\ref{70z})$ and $(\ref{70'z})$ with making $y=\frac{2}{2-\alpha}\mu x^{\frac{2-\alpha}{2}}$ and $X=\left( \frac{2-\alpha}{2\mu}s\right) ^{\frac{2}{2-\alpha}}$, we get
\begin{equation}
\matrix{v(x)&=& A\theta_+(x)+B\theta_-(x)\hfill \cr
&+&\frac{\pi}{2\sin \nu_\alpha \pi}
\left(\frac{2}{2-\alpha}\right)\int_{0}^{x}\imath f_1(X)(\theta_+(X)
\theta_-(x)-\theta_+(x)\theta_-(X))dX,\hfill \cr}
\label{eqtt}
\end{equation}
where $\theta_+$ and $\theta_-$ are defined by $(\ref{theta})$,
then
\begin{equation}
\matrix{v_x(x)&=& A\theta'_+(x)+B\theta'_-(x)\hfill \cr
&+&\frac{\pi}{2\sin \nu_\alpha \pi}
\left(\frac{2}{2-\alpha}\right)\int_{0}^{x}\imath f_1(X)(\theta_+(X)
\theta'_-(x)-\theta'_+(x)\theta_-(X))dX.\hfill \cr}
\label{eqtt1}
\end{equation}
To reformulate the conditions $(\ref{cond})_2$ and $(\ref{cond22})$ we use the expressions of $v$ and $v_x$.\\

The first boundary condition (\ref{cond22}) become
$$
A(1-\alpha)d^+-\imath\rho(i\lambda+\eta)^{\tilde{\alpha}-1} Bd^-=i\zeta\Int_{-\infty}^{+\infty}\Frac{\mu(\xi)f_2(\xi)}{i\lambda+\xi^{2}+\eta}\, d\xi.,
$$
where $d^+$ and $d^-$ are defined by $(\ref{dd})$.

The second condition $v_x(1)=0$ become
$$
A\theta'_+(1)+B\theta'_-(1)=-\frac{\pi}{2\sin \nu_\alpha \pi}
\left(\frac{2}{2-\alpha}\right)\int_{0}^{1}\imath f_1(X)(\theta_+(X)
\theta'_-(1)-\theta'_+(1)\theta_-(X))dX.
$$
In order to get the expressions of $\theta_+'(1)$ and $\theta_-'(1)$, we derive $\theta_+$ and $\theta_-$ respectively
and we use the following relation
\begin{equation}\label{45}
xJ'_{\nu_{\alpha}}=\nu_{\alpha}J_{\nu_{\alpha}}(x)-xJ_{\nu_{\alpha}+1}(x),
\end{equation}
we deduce that
\begin{equation}\label{theta+'}
\theta_+'(1)=(1-\alpha)J_{\nu_{\alpha}}\left(\frac{2\mu}{2-\alpha}\right)-\mu J_{\nu_{\alpha}+1}\left(\frac{2\mu}{2-\alpha}\right)
\end{equation}
and
\begin{equation}\label{theta-'}
\theta_-'(1)=-\mu J_{-\nu_{\alpha}+1}\left(\frac{2\mu}{2-\alpha}\right).
\end{equation}
Therefore, we get the following linear system in $A$ and $B$
\begin{equation}\label{77z}
\left( \begin{array}{cc}
(1-\alpha)d^+ &
-\imath\rho(i\lambda+\eta)^{\tilde{\alpha}-1} d^- \\\\
\theta'_+(1) & \theta'_-(1)
\end{array}\right)
\left( \begin{array}{c}
A \\\\ B
\end{array} \right) =
\left( \begin{array}{c}
C \\\\ \tilde{C}
\end{array} \right) ,
\end{equation}
As $D\not= 0$ for all $\lambda\not=0$, then $A$ and $B$ are uniquely determined by (\ref{77z}).
Now, we prove that $(v, \phi)\in D({\cal A})$.

\noindent
First $V\in H_{\alpha}^{2}$. Indeed,
from (\ref{theta}), we have
\begin{equation}
\left\{\matrix{x^{\alpha/2}\theta'_{+}(x)&=&(1-\alpha) x^{-1/2}J_{\nu_\alpha}\left(\frac{2}{2-\alpha}\mu x^{\frac{2-\alpha}{2}}\right)
-\mu x^{\frac{1-\alpha}{2}}J_{1+\nu_\alpha}\left(\frac{2}{2-\alpha}\mu x^{\frac{2-\alpha}{2}}\right),\hfill \cr
x^{\alpha/2}\theta'_{-}(x)&=&
-\mu x^{\frac{1-\alpha}{2}}J_{1-\nu_\alpha}\left(\frac{2}{2-\alpha}\mu x^{\frac{2-\alpha}{2}}\right).\hfill \cr
(x^{\alpha}\theta'_{+})'(x)&=&-(3-2\alpha)\mu x^{-1/2}J_{1+\nu_\alpha}\left(\frac{2}{2-\alpha}\mu x^{\frac{2-\alpha}{2}}\right)
+\mu^2 x^{\frac{1-\alpha}{2}}J_{2+\nu_\alpha}\left(\frac{2}{2-\alpha}\mu x^{\frac{2-\alpha}{2}}\right),\hfill \cr
(x^{\alpha}\theta'_{-})'(x)&=&-\mu x^{-1/2}J_{1-\nu_\alpha}\left(\frac{2}{2-\alpha}\mu x^{\frac{2-\alpha}{2}}\right)
+\mu^2 x^{\frac{1-\alpha}{2}}J_{2-\nu_\alpha}\left(\frac{2}{2-\alpha}\mu x^{\frac{2-\alpha}{2}}\right),\hfill \cr
}\right.
\label{eqtt2}
\end{equation}
In the following Lemma we will give some technical inequalities which will be useful
for showing our results.
\begin{lemma}\label{l5.3}
We have
\begin{equation}\label{80}
\|\theta_+\|_{L^2(0,1)},\|\theta_-\|_{L^2(0,1)}\leq\frac{c}{\sqrt{\mu}}.
\end{equation}
\begin{equation}
\left\|x^{-\frac{1}{2}}
J_{\nu_{\alpha}}\left(\Frac{2}{2-\alpha}\mu x^{\frac{2-\alpha}{2}}\right)\right\|_{L^{2}(0, 1)}, \left\|x^{-\frac{1}{2}}
J_{-\nu_{\alpha}}\left(\Frac{2}{2-\alpha}\mu x^{\frac{2-\alpha}{2}}\right)\right\|_{L^{2}(0, 1)}\leq c\sqrt{|\mu|}.
\label{hhgg5}
\end{equation}
\end{lemma}
The proof of Lemma \ref{l5.3} will be given in Appendix $A$.\\

Now, using (\ref{eqtt}), (\ref{eqtt1}) and (\ref{eqtt2}), it easy to see that $v\in H_{\alpha}^{2}(0, 1)$.
Moreover $\phi, \xi\phi \in L^2(\R)$.\\

\noindent
\textbf{Case $2$}: $\lambda=0$.
We can obtain the result using the Lax-Milgram Theorem.

\hfill$\diamondsuit$\\

According to the Lemmas $\ref{l5.1}$, $\ref{l5.2}$ and $\ref{5.1}$ the $C_0$-semigroup
$e^{t{\mathcal{A}}}$ is strongly stable in ${\mathcal{H}}$.\\

Next, in order to prove an polynomial decay rate we will use the following theorem.

\begin{lemma}\label{5.2}{\bf\cite{BT}}
Let $S(t)$ be a bounded $C_0$-semigroup on a Hilbert space ${\cal X}$ with generator ${\cal A}$. If
$$
\imath \R\subset\rho({\cal A}) \hbox{ and } \
\overline{\lim_{|\beta|\rightarrow\infty}}\frac{1}{\beta^l}
\|(\imath\beta I-{{\cal A}})^{-1}\|_{{\cal L}({\cal X})}<\infty
$$
for some $l$, then there exist $c$ such that
$$
\|e^{{\cal A}t}V_0\|^2\leq \frac{c}{t^{\frac{2}{l}}}\|V_0\|^2_{D({\cal A})}.
$$
\end{lemma}

Our main result is the following.
\begin{theorem}\label{5.4}
If $\eta\not=0$, then the global solution of the problem $(\ref{P})$ has the following energy decay property
$$
E(t)=\|S_{\cal A}(t)V_0\|_{\cal H}^2\leq\left\{\matrix{
\Frac{c}{ t^{\frac{2}{\nu_\gamma-\tilde{\alpha}+\frac{1}{2}}}}\|V_0\|_{D({\cal A})}^2 &\hbox{ if }{\tilde\alpha}<\frac{4-3\alpha}{2(2-\alpha)},\hfill \cr
c e^{-\omega t}\|V_0\|_{{\cal H}}^2&\hbox{ if }{\tilde\alpha}\geq\frac{4-\alpha}{2(2-3\alpha)},\hfill \cr}\right.
$$
where $\nu_\alpha=\frac{1-\alpha}{2-\alpha}$.
Moreover, the rate of energy decay is optimal for general initial data in $D({\cal A})$.
\end{theorem}
{\bf Proof.}
We need to estimate $\|V\|_{{\cal H}}$, where $V$ is a solution of the resolvent equation given by
$$
(\imath\lambda-{\mathcal{A}})V=F,
$$
where $\lambda\in \R$ and $F\in \mathcal{H}$.

Throughout this proof we use the notation introduced in the proof of
Lemma \ref{l5.2}. Inverting the matrix (\ref{77z}) we obtain
\begin{equation}\label{77}
\left( \begin{array}{c}
A \\\\ B
\end{array} \right)=
\left( \begin{array}{cc}
\theta'_-(1) &
\imath\rho(i\lambda+\eta)^{\tilde{\alpha}-1} d^- \\\\
-\theta'_+(1) &(1-\alpha)d^+
\end{array}\right)
\left( \begin{array}{c}
C \\\\ \tilde{C}
\end{array} \right) ,
\end{equation}
where
$$
\theta_+'(1)=(1-\alpha)J_{\nu_{\alpha}}\left(\frac{2\mu}{2-\alpha}\right)-\mu J_{\nu_{\alpha}+1}\left(\frac{2\mu}{2-\alpha}\right),
$$
$$
\theta_-'(1)=-\mu J_{-\nu_{\alpha}+1}\left(\frac{2\mu}{2-\alpha}\right),
$$
and
$$
C=i\zeta\Int_{-\infty}^{+\infty}\Frac{\mu(\xi)f_2(\xi)}{i\lambda+\xi^{2}+\eta}\, d\xi,
$$
and
$$
\tilde{C}=-\frac{\pi}{2\sin \nu_\alpha \pi}
\left(\frac{2}{2-\alpha}\right)\int_{0}^{1}\imath f_1(X)(\theta_+(X)
\theta'_-(1)-\theta'_+(1)\theta_-(X))dX.
$$
Let $D$ the determinant of the linear system $(\ref{77z})$, so using $(\ref{84})$ we get
\begin{eqnarray*}
D&=&-(1-\alpha)c^+_{\nu_{\alpha,0}}\left(\frac{2}{2-\alpha}\right)^{\nu_{\alpha}}\left(\frac{2-\alpha}{\pi}\right)^{\frac{1}{2}}
\mu^{\nu_{\alpha}+\frac{1}{2}}\cos\left(\frac{2\mu}{2-\alpha}-(1-\nu_{\alpha})\frac{\pi}{2}-\frac{\pi}{4} \right)\\
&+&\imath\rho (1-\alpha)(i\lambda+\eta)^{\tilde{\alpha}-1} c^-_{\nu_{\alpha,0}}\left(\frac{2}{2-\alpha}\right)^{-\nu_{\alpha}}\left(\frac{2-\alpha}{\pi}\right)^{\frac{1}{2}}
\mu^{-\nu_{\alpha}-\frac{1}{2}}\cos\left(\frac{2\mu}{2-\alpha}+\nu_{\alpha}\frac{\pi}{2}-\frac{\pi}{4} \right)\\
&-&\imath\rho (i\lambda+\eta)^{\tilde{\alpha}-1} c^-_{\nu_{\alpha,0}}\left(\frac{2}{2-\alpha}\right)^{-\nu_{\alpha}}\left(\frac{2-\alpha}{\pi}\right)^{\frac{1}{2}}
\mu^{-\nu_{\alpha}+\frac{1}{2}}\cos\left(\frac{2\mu}{2-\alpha}-(\nu_{\alpha}+1)\frac{\pi}{2}-\frac{\pi}{4} \right)+O\left( \frac{1}{|\mu|}\right).
\end{eqnarray*}
It is clear that
\begin{equation}\label{78}
|D|\geq c|\mu|^{2\tilde{\alpha}-\nu_{\alpha}-\frac{3}{2}},\ \ \hbox{for large } |\mu|.
\end{equation}
Indeed, suppose (\ref{78}) was wrong. Then $\exists\mu_n$ such that $|\mu_n|\rightarrow \infty$ with
\begin{equation}\label{78sd}
|D||\mu_n|^{-2\tilde{\alpha}+\nu_{\alpha}+\frac{3}{2}}\rightarrow 0 \hbox{ as }n\rightarrow +\infty.
\end{equation}
By $\Re D$,
$$
|\mu_n|^{-2\tilde{\alpha}+\nu_{\gamma}+2}\cos\left(\frac{2}{2-\gamma}\mu_n-(1-\nu_{\gamma})\frac{\pi}{2}-\frac{\pi}{4}\right) \rightarrow 0 \hbox{ as } n\rightarrow +\infty.
$$
By $\Im D$,
$$
\cos\left(\frac{2}{2-\gamma}\mu_n-(1+\nu_{\gamma})\frac{\pi}{2}-\frac{\pi}{4}\right) \rightarrow 0 \hbox{ as } n\rightarrow +\infty.
$$
This is impossible. Indeed, $\exists k_n\in \Z$ with $|k_n|\rightarrow +\infty\ n\rightarrow +\infty$ such that
$$
\frac{2}{2-\gamma}\mu_n-(1+\nu_{\gamma})\frac{\pi}{2}-\frac{\pi}{4}=(k_n+\Frac{1}{2})\pi+o(1).
$$
Then
$$
\left|\cos\left(\frac{2}{2-\gamma}\mu_n-(1-\nu_{\gamma})\frac{\pi}{2}-\frac{\pi}{4}\right)\right|
\rightarrow \sin\nu_{\gamma}\pi \hbox{ as }n\rightarrow +\infty.
$$
According the the linear system $(\ref{77})$, we have
$$
|A|=\left| \frac{\theta_-'(1)C+\imath\rho(i\lambda+\eta)^{\tilde{\alpha}-1} \tilde{C}d^-}{D} 
\right|
$$
and
$$
|B|=\left| \frac{-\theta_+'(1)C+(1-\alpha)\tilde{C}d^+}{D} 
\right|.
$$
In order to estimate $\tilde{C}$, we use Lemma \ref{l5.3}, where in which we consider $\mu>0$.

From Lemma \ref{l5.3} and the asymptotic formula $(\ref{84})$ for large $\mu$, we deduce that
$$
\matrix{\|\theta_+\|_{L^2(0,1)}, \|\theta_-\|_{L^2(0,1)}&\leq &\frac{c}{\sqrt{|\mu|}}, \hfill \cr
|\theta'_+(1)|, |\theta'_-(1)|&\leq & c \sqrt{|\mu|}.\hfill \cr}
$$

\hfill$\diamondsuit$\\

Then, using Cauchy-Schwartz inequality, the expressions of $\theta_+'$ and $\theta_-'$, we get
\begin{equation}\label{qqa}
\matrix{|C| &\leq & \zeta \Int_{-\infty}^{\infty}\Frac{\mu(\xi)^2}{|i\lambda+\xi^2+\eta|^2}\, d\xi\|f_2\|_{L^2(0,1)}\hfill \cr
&\leq & 2\zeta \Int_{-\infty}^{\infty}\Frac{\mu(\xi)^2}{(|\lambda|+\xi^2+\eta)^2}\, d\xi\|f_2\|_{L^2(0,1)}\hfill \cr
&\leq & c|\mu|^{\tilde{\alpha}-2}\|f_2\|_{L^2(0,1)}.\hfill \cr}
\end{equation}
\begin{equation}\label{C}
|\tilde{C}|\leq c\|f_1\|_{L^2(0,1)}.
\end{equation}
Using $(\ref{C})$, $(\ref{dd})$, and $(\ref{78})$, we deduce that
$$
|A|\leq c|\mu|^{-\frac{1}{2}}+c'|\mu|^{\nu_{\alpha}-\tilde{\alpha}}\
$$
and
$$
\matrix{|B|&\leq & c |\mu|^{\nu_{\alpha}-\tilde{\alpha}}+c'|\mu|^{2\nu_{\alpha}-2\tilde{\alpha}+\frac{3}{2}}\hfill \cr
&\leq & c'|\mu|^{2\nu_{\alpha}-2\tilde{\alpha}+\frac{3}{2}}.\hfill \cr}
$$
From the expression of $v$, we deduce that
$$
\|v\|_{L^2(0,1)}\leq c'|\mu|^{2\nu_{\alpha}-2\tilde{\alpha}+1} \|F\|_{\cal H}.
$$
Since $\eta> 0$, we have (see (\ref{zz18}))
\begin{equation}
\|\phi\|_{L^{2}(-\infty, \infty)}^2\leq \Int_{-\infty}^{+\infty}(\xi^{2}+\eta)|\phi(\xi)|^{2}\, d\xi\leq c\|V\|_{{\cal H}}\|F\|_{{\cal H}}.
\label{aaqq1}
\end{equation}
Thus, we conclude that
$$
\|v\|_{L^2(0,1)}^2+\|\phi\|_{L^{2}(-\infty, \infty)}^2\leq c'|\mu|^{2(2\nu_{\alpha}-2\tilde{\alpha}+1)} \|F\|_{\cal H}^2
+c\|V\|_{{\cal H}}\|F\|_{{\cal H}}.
$$
Hence
$$
\|V\|_{{\cal H}}^2\leq c'|\mu|^{2(2\nu_{\alpha}-2\tilde{\alpha}+1)} \|F\|_{\cal H}^2+c''\|F\|_{{\cal H}}^2.
$$
So
$$
\|(\imath\lambda-{\mathcal{A}})^{-1}\|_{\mathcal{H}}\leq c \left\{\matrix{|\lambda|^{\nu_{\alpha}-\tilde{\alpha}+\frac{1}{2}} \ \
\hfill &\hbox{ as } \lambda \rightarrow \infty \hfill &\hbox{ if } \nu_{\alpha}-\tilde{\alpha}+\frac{1}{2}> 0 \hfill \cr
C \hfill &\hbox{ as } \lambda \rightarrow \infty \hfill& \hbox{ if } \nu_{\alpha}-\tilde{\alpha}+\frac{1}{2}\leq 0 \hfill \cr}\right.
$$

The conclusion follows by applying Lemma $\ref{5.2}$.
\begin{remark}
It possible to obtain a charp estimate of the resolvent in the case $\nu_{\alpha}-\tilde{\alpha}+\frac{1}{2}\leq 0$.
Indded instead (\ref{aaqq1}), we use (\ref{58}). We have
$$
\phi(\xi)=\Frac{v(0)\mu(\xi)}{1\lambda+\xi^2+\eta}+\Frac{f_2(\xi)}{1\lambda+\xi^2+\eta}.
$$
Then
$$
\|\phi\|_{L^{2}(0, 1)}^2\leq C |v(0)|^2 |\lambda|^{\tilde{\alpha}-2}+\Frac{c}{|\lambda|^2}\|f_2\|_{L^2({\R})}^2
$$
From (\ref{eqtt}), we have
$$
\matrix{|v(0)|^2&\leq & c |B|^2 |\mu|^{-2\nu_{\alpha}}\hfill \cr
&\leq & c |\lambda|^{\nu_{\alpha}-2\tilde{\alpha}+3/2}\|f_1\|_{L^2(0, 1)}^2\hfill \cr}
$$
Hence
$$
\matrix{\|\phi\|_{L^{2}(0, 1)}^2&\leq & c |\lambda|^{\nu_{\alpha}-\tilde{\alpha}-1/2}\|f_1\|_{L^2(0, 1)}^2
+\Frac{c}{|\lambda|^2}\|f_2\|_{L^2({\R})}^2\hfill \cr
&\leq & c |\lambda|^{\nu_{\alpha}-\tilde{\alpha}-1/2}\|F\|_{\cal H}^2.\hfill \cr}
$$
Finally, we conclude
$$
\|V\|_{{\cal H}}^2\leq c'|\mu|^{2(2\nu_{\alpha}-2\tilde{\alpha}+1)} \|F\|_{\cal H}^2\ \hbox{ if } \nu_{\alpha}-\tilde{\alpha}+\frac{1}{2}\leq 0
$$
and so
$$
\|(\imath\lambda-{\mathcal{A}})^{-1}\|_{\mathcal{H}}\leq c |\lambda|^{\nu_{\alpha}-\tilde{\alpha}+\frac{1}{2}}\ \hbox{ if } \nu_{\alpha}-\tilde{\alpha}+\frac{1}{2}\leq 0.
$$
\end{remark}

\hfill$\diamondsuit$\\

\section{Optimality of energy decay}

In this section, we will study the lack of exponential decay of solution of the system $(\ref{10})$. For this purpose
we will use the following theorem.
\begin{lemma}\label{4.1}{\bf\cite{P}}
Let $S(t)$ be a $C_0$-semigroup of contractions on Hilbert space $\mathcal{X}$ with generator ${\mathcal{A}}$. Then $S(t)$ is exponentially stable if and only if
$$
\rho({\mathcal{A}})\supseteq\{\imath\beta: \beta \in \R\}\equiv\imath\R
$$
and
$$
\overline{\lim_{|\beta|\rightarrow\infty}}
\|(\imath\beta I-{\mathcal{A}})^{-1}\|_{\mathcal{L(X)}}<\infty.
$$
\end{lemma}

Our main result is the following.
\begin{theorem}\label{4.2}
The semigroup generated by the operator ${\mathcal{A}}$ is not exponentially stable for\\ $\tilde{\alpha}< \frac{4-3\alpha}{2(2-\alpha)}$.
\end{theorem}
{\bf Proof.}
We aim to show that an infinite number of eigenvalues of ${\mathcal{A}}$ approach the imaginary axis
which prevents the system $(\ref{P})$ from being exponentially stable. Let $\lambda$ be an eigenvalue of ${\mathcal{A}}$ with associated eigenvector $v$. Then the equation ${\mathcal{A}}v=\lambda v$ is equivalent to
$$
\imath\lambda v-(x^\alpha v_x)_x=0
$$
together with the conditions
$$
\quad \left \{
\begin{array}{ll}
(x^\alpha v_x)(0)=\imath\rho (\lambda+\eta)^{\tilde{\alpha}-1}v(0),\\
v_x(1)=0,
\end{array}
\right.
$$
so we get the following system
\begin{equation}\label{40}
\quad \left \{
\begin{array}{ll}
\gamma^2v-(x^\alpha v_x)_x=0,\\
(x^\alpha v_x)(0)=\imath\rho (\lambda+\eta)^{\tilde{\alpha}-1}v(0),\\
v_x(1)=0,
\end{array}
\right.
\end{equation}
with $\gamma^2=\imath\lambda$.\\

Suppose that $v$ is a solution of $(\ref{40})_1$, then the function $\Psi$ defined by
$$ 
v(x)=x^{\frac{1-\alpha}{2}}\Psi\left(\frac{2}{2-\alpha} \imath\gamma
x^{\frac{2-\alpha}{2}} \right)
$$ 
is a solution of the following equation
\begin{equation}\label{41}
y^2\Psi''(y)+y\Psi'(y)+\left(y^2-\left( \frac{\alpha-1}{2-\alpha}\right) ^2 \right) \Psi(y)=0.
\end{equation}
We have
\begin{equation}\label{42}
v(x)=c_+\widetilde{\theta}_+(x)+c_-\widetilde{\theta}_-(x)
\end{equation}
where
$$
\widetilde{\theta}_+(x)=x^{\frac{1-\alpha}{2}}J_{\nu_{\alpha}}\left(\frac{2}{2-\alpha}\imath\gamma x^{\frac{2-\alpha}{2}}\right)
\quad \hbox{ and  } \ \ \
\widetilde{\theta}_-(x)=x^{\frac{1-\alpha}{2}}J_{-\nu_{\alpha}}\left(\frac{2}{2-\alpha}\imath\gamma x^{\frac{2-\alpha}{2}}\right).
$$
Therefore the boundary conditions can be written as the following system
\begin{equation}\label{46}
\widetilde{M}(\gamma)C(\gamma)=\left( \begin{array}{cc}
(1-\alpha)\widetilde{d}^+ &
-\imath\rho (\lambda+\eta)^{\tilde{\alpha}-1}\widetilde{d}^- \\\\
\widetilde{\theta}'_+(1) & \widetilde{\theta}'_-(1)
\end{array}\right)
\left( \begin{array}{c}
c_+ \\\\ c_-
\end{array} \right) =
\left( \begin{array}{c}
0 \\\\ 0
\end{array} \right) ,
\end{equation}
where
\begin{equation}\label{ddtilde}
\widetilde{d}^+=c^+_{\nu_{\alpha,0}}\left(\frac{2}{2-\alpha}\imath\gamma\right)^{\nu_{\alpha}},\ \
\widetilde{d}^-=c^-_{\nu_{\alpha,0}}\left(\frac{2}{2-\alpha}\imath\gamma\right)^{{-\nu}_{\alpha}},
\end{equation}
\begin{equation}\label{theta'+tilde}
\widetilde{\theta}_+'(1)=(1-\alpha)J_{\nu_{\alpha}}\left(\frac{2\gamma}{2-\alpha}\imath\right)-\imath\gamma J_{\nu_{\alpha}+1}\left(\frac{2\gamma}{2-\alpha}\imath\right)
\end{equation}
and
\begin{equation}\label{theta'-tilde}
\widetilde{\theta}'_-(1)=-\imath\gamma J_{-\nu_{\alpha}+1}\left(\frac{2\gamma}{2-\alpha}\imath\right).
\end{equation}
Then, a non-trivial solution $v$ exists if and only if the determinant of $\widetilde{M}(\gamma)$ vanishes.
Set $f(\gamma)=\det\widetilde{M}(\gamma)$, thus the characteristic equation is $f(\gamma)=0$.

Our purpose is to prove, thanks to Rouch\'{e}'s Theorem, that there is a subsequence of eigenvalues for which their
real part tends to $0$.

In the sequel, since ${\mathcal{A}}$ is dissipative, we study the asymptotic behavior of the large eigenvalues $\lambda$
of ${\mathcal{A}}$ in the strip $-\alpha_0\leq\Re(\lambda)\leq0$, for some $\alpha_0>0$ large enough and for such
$\lambda$, we remark that $\theta_+$ and $\theta_-$ remain bounded.
\begin{lemma}\label{l4.1}
There exists $N \in \N$ such that
$$
\{\lambda_k\}_{{k\in\Z}^*, \ |k|\geq N} \subset \sigma({\mathcal{A}}),
$$
where

\noindent
$\bullet$ If $\tilde{\alpha}=1$, then
\begin{eqnarray*}
\lambda_k=i \left[C_0^2(k\pi)^2+2 C_0 C_1 k\pi^2\right]
+2\Frac{C_0 C_2 \sin\nu_{\alpha} \pi}{(k\pi)^{2\nu_{\alpha}-1}}+O(1)
\end{eqnarray*}
where
$$
C_0=-\frac{2-\alpha}{2},\quad C_1=-\frac{2-\alpha}{2}\left(-\frac{\nu_{\alpha}}{2}+\frac{5}{4}\right),
$$
$$
C_2=\Frac{\rho(2-\alpha)}{2(1-\alpha)}\frac{c^-_{\nu_{\alpha,0}}}{c^+_{\nu_{\alpha,0}}}.
$$

\noindent
$\bullet$ If $\tilde{\alpha}>\Frac{4-3\alpha}{2(2-\alpha)}$, then
\begin{eqnarray*}
\lambda_k=i \left[C_0^2(k\pi)^2+2 C_0 C_1 k\pi^2\right]
-2i\Frac{C_0 C_2 (-i)^{3\tilde{\alpha}} \sin\nu_{\alpha} \pi}{(k\pi)^{2\nu_{\alpha}-2\tilde{\alpha}-1}}+O(1)
\end{eqnarray*}
where
$$
C_0=-\frac{2-\alpha}{2},\quad C_1=-\frac{2-\alpha}{2}\left(-\frac{\nu_{\alpha}}{2}+\frac{5}{4}\right),
$$
$$
C_2=\Frac{\rho(2-\alpha)}{2(1-\alpha)}\left(\Frac{2}{2-\alpha}\right)^{2-2\tilde{\alpha}}\frac{c^-_{\nu_{\alpha,0}}}{c^+_{\nu_{\alpha,0}}}.
$$
\noindent
$\bullet$ If $\nu_{\alpha} <\tilde{\alpha}<\Frac{4-3\alpha}{2(2-\alpha)}$, then
\begin{eqnarray*}
\lambda_k=&=&i \left[C_0^2(k\pi)^2+C_1^2\pi^2+2 C_0 C_1 k\pi^2+2{C_0 C_3}\right]
-2i\Frac{C_0 C_2(-i)^{3\tilde{\alpha}}\sin\nu_{\alpha} \pi}{(k\pi)^{2\nu_{\alpha}-2\tilde{\alpha}+1}}+O\left(\frac{1}{k}\right),
\end{eqnarray*}
where
$$
C_0=-\frac{2-\alpha}{2},\quad C_1=-\frac{2-\alpha}{2}\left(-\frac{\nu_{\alpha}}{2}+\frac{5}{4}\right),
$$
$$
C_2=\Frac{\rho(2-\alpha)}{2(1-\alpha)}\left(\Frac{2}{2-\alpha}\right)^{2-2\tilde{\alpha}}\frac{c^-_{\nu_{\alpha,0}}}{c^+_{\nu_{\alpha,0}}},\quad
C_3=\Frac{(2-\alpha)}{4}{(\frac{1}{2}-\nu_{\alpha})(\frac{3}{2}-\nu_{\alpha})}.
$$
\noindent
$\bullet$ If $\tilde{\alpha}<\nu_{\alpha}$, then
\begin{eqnarray*}
\lambda_k=i \left[C_0^2(k\pi)^2+C_1^2\pi^2+2 C_0 C_1 k\pi^2+2{C_0 C_2}+2\frac{C_0 C_3}{k}+2\frac{C_0 C_1}{k}\right]
-2i\Frac{C_0 C_4(-i)^{3\tilde{\alpha}}\sin\nu_{\alpha} \pi}{(k\pi)^{2\nu_{\alpha}-2\tilde{\alpha}+1}}
+O\left(\frac{1}{k^{2}}\right)
\end{eqnarray*}
where
$$
C_0=-\frac{2-\alpha}{2},\quad C_1=-\frac{2-\alpha}{2}\left(-\frac{\nu_{\alpha}}{2}+\frac{5}{4}\right),\quad
C_2=\Frac{(2-\alpha)}{4}(\frac{1}{2}-\nu_{\alpha})(\frac{3}{2}-\nu_{\alpha}),
$$
and
$$
C_3=-\frac{2-\alpha}{2}m (-\frac{\nu_{\alpha}}{2}+\frac{5}{4}),\quad
C_4=\Frac{\rho(2-\alpha)}{2(1-\alpha)}\left(\Frac{2}{2-\alpha}\right)^{2-2\tilde{\alpha}}\frac{c^-_{\nu_{\alpha,0}}}{c^+_{\nu_{\alpha,0}}}.
$$
and
$$
\lambda_k=\overline{\lambda_{-k}}, \hbox{ if } k \leq -N.
$$
\end{lemma}
{\bf Proof.} $\bullet\ \ \tilde{\alpha}=1$.
We aim to solve the equation
$$
f(\gamma)=-\imath\gamma(1-\alpha)\widetilde{d}^+
J_{1-\nu_\alpha}\left(\frac{2\gamma}{2-\alpha}\imath\right)+
\imath\rho(1-\alpha)\widetilde{d}^-
J_{\nu_\alpha}\left(\frac{2\gamma}{2-\alpha}\imath\right)
+\rho{\widetilde{d}}^-\gamma
J_{1+\nu_\alpha}\left(\frac{2\gamma}{2-\alpha}\imath\right)=0.
$$
We will use the following classical development (see {\bf\cite{Le}}), for all $\delta>0$ and when $|\arg{z}|<\pi-\delta$:
\begin{equation}
\matrix{J_{\nu}(z)=\left(\Frac{2}{\pi z}\right)^{1/2}\left[\cos\left(z-\nu\frac{\pi}{2}-\frac{\pi}{4}\right)
-\Frac{(\nu-\Frac{1}{2})(\nu+\Frac{1}{2})}{2}\Frac{\sin\left(z-\nu\frac{\pi}{2}-\frac{\pi}{4}\right)}{z}\right.\hfill & \cr
\left.-\Frac{(\nu-\Frac{1}{2})(\nu+\Frac{1}{2})(\nu-\Frac{3}{2})(\nu+\Frac{3}{2})}{8}\Frac{\cos\left(z-\nu\frac{\pi}{2}-\frac{\pi}{4}\right)}{z^2}
+O\left(\frac{1}{|z|^3}\right)\right]. & \cr}
\label{84}
\end{equation}
We get
$$
f(\gamma)=-\imath(1-\alpha)c^+_{\nu_{\alpha,0}}
\gamma^{1+\nu_{\alpha}}\left( \frac{2}{2-\alpha}\imath\right)^{\nu_{\alpha}}
\left( \frac{2}{\pi z}\right)^\frac{1}{2}
\frac{e^{-\imath(z+(\nu_{\alpha}-1)\frac{\pi}{2}-\frac{\pi}{4})}}{2}\widetilde{f}(\gamma),
$$
where $$z=\frac{2\gamma}{2-\alpha}\imath$$
and
\begin{eqnarray*}
\widetilde{f}(\gamma)&=& 1+e^{2\imath(z+(\nu_{\alpha}-1)\frac{\pi}{2}-\frac{\pi}{4})}-
\frac{\rho}{\imath}\frac{1}{1-\alpha}\frac{c^-_{\nu_{\alpha,0}}}{c^+_{\nu_{\alpha,0}}}\left( \frac{2}{2-\alpha}\imath\right) ^{-2\nu_{\alpha}}
\frac{e^{\imath(2z-3\frac{\pi}{2})}+e^{\imath\nu_{\alpha}\pi}}{\gamma^{2\nu_{\alpha}}}\\\\
&-&\Frac{(\frac{1}{2}-\nu_{\alpha})(\frac{3}{2}-\nu_{\alpha})}{2i}\left( \frac{2}{2-\alpha}\imath\right)^{-1}
\Frac{e^{2\imath(z+(\nu_{\alpha}-1)\frac{\pi}{2}-\frac{\pi}{4})}-1}{\gamma}\\ \\
&-&\rho\frac{1}{1-\alpha}\frac{c^-_{\nu_{\alpha,0}}}{c^+_{\nu_{\alpha,0}}}
\left( \frac{2}{2-\alpha}\imath\right) ^{-2\nu_{\alpha}-1}\Frac{(\nu_{\alpha}+\frac{1}{2})(\nu_{\alpha}+\frac{3}{2})}{2}
\frac{e^{\imath(2z-3\frac{\pi}{2})}-e^{\imath\nu_{\alpha}\pi}}{\gamma^{1+2\nu_{\alpha}}}\\ \\
&-&\rho \frac{c^-_{\nu_{\alpha,0}}}{c^+_{\nu_{\alpha,0}}}
\left( \frac{2}{2-\alpha}\imath\right) ^{-2\nu_{\alpha}}\frac{e^{\imath(2z-\pi)}+e^{\imath(\nu_{\alpha}\pi-\frac{\pi}{2})}}{\gamma^{1+2\nu_{\alpha}}}+O\left(\frac{1}{\gamma^2} \right)\\\\
&=&f_0(\gamma)+\frac{f_1(\gamma)}{\gamma^{2\nu_{\alpha}}}+\frac{f_2(\gamma)}{\gamma}+\frac{f_3(\gamma)}{\gamma^{1+2\nu_{\alpha}}}+O\left(\frac{1}{\gamma^2} \right)
\end{eqnarray*}
with
$$
\matrix{f_0(\gamma)&=&1+e^{2\imath(z+(\nu_{\alpha}-1)\frac{\pi}{2}-\frac{\pi}{4})}, \hfill \cr
f_1(\gamma)&=&-\frac{\rho}{\imath}\frac{1}{1-\alpha}\frac{c^-_{\nu_{\alpha,0}}}{c^+_{\nu_{\alpha,0}}}\left( \frac{2}{2-\alpha}\imath\right) ^{-2\nu_{\alpha}}
(e^{\imath(2z-3\frac{\pi}{2})}+e^{\imath\nu_{\alpha}\pi}),\hfill \cr
f_2(\gamma)&=&-\Frac{(\frac{1}{2}-\nu_{\alpha})(\frac{3}{2}-\nu_{\alpha})}{2i}\left( \frac{2}{2-\alpha}\imath\right)^{-1}
(e^{2\imath(z+(\nu_{\alpha}-1)\frac{\pi}{2}-\frac{\pi}{4})}-1).\hfill \cr
f_3(\gamma)&=&-\rho \frac{c^-_{\nu_{\alpha,0}}}{c^+_{\nu_{\alpha,0}}}
\left( \frac{2}{2-\alpha}\imath\right) ^{-2\nu_{\alpha}}(e^{\imath(2z-\pi)}+e^{\imath(\nu_{\alpha}\pi-\frac{\pi}{2})})\hfill  \cr
&&\ \ \ -\rho\frac{1}{1-\alpha}\frac{c^-_{\nu_{\alpha,0}}}{c^+_{\nu_{\alpha,0}}}
\left( \frac{2}{2-\alpha}\imath\right) ^{-2\nu_{\alpha}-1}\Frac{(\nu_{\alpha}+\frac{1}{2})(\nu_{\alpha}+\frac{3}{2})}{2}
(e^{\imath(2z-3\frac{\pi}{2})}-e^{\imath\nu_{\alpha}\pi}),\hfill & \cr}
$$
Note that $f_0, f_1, f_2$ and $f_3$ remain bounded in the strip $-\alpha_0\leq\Re(\lambda)\leq0$.\\

We search the roots of $f_0$,
$$
f_0(\gamma)=0 \ \Leftrightarrow \ 1+e^{2\imath(z+(\nu_{\alpha}-1)\frac{\pi}{2}-\frac{\pi}{4})}=0,
$$
so, $f_0$ has the following roots
$$
\gamma_k^0=-\frac{2-\alpha}{2}\imath\left(k-\frac{\nu_{\alpha}}{2}+\frac{5}{4}\right)\pi, \ \ k \in \Z.
$$
Let $B_k(\gamma_k^0,r_k)$ be the ball of centrum $\gamma_k^0$ and
radius $r_k=\frac{1}{k^{\nu_{\alpha}}}$, then if $\gamma\in \partial B_k$, we have $\gamma=\gamma_k^0+r_ke^{\imath\theta}, \ \theta\in [0,2\pi] $, then we have
$$
f_0(\gamma)=\frac{4}{2-\alpha}r_ke^{\imath\theta}+o(r_k^2).
$$
Hence, there exists a positive constant c such that, for all $\gamma\in \partial B_k$
$$
|f_0(\gamma)|\geq cr_k=\frac{c}{k^{\nu_{\alpha}}}.
$$

From the expression of $\widetilde{f}$, we conclude that
$$
|\widetilde{f}(\gamma)-f_0(\gamma)|=O\left(\frac{1}{\gamma^{2\nu_{\alpha}}} \right) =O\left(\frac{1}{k^{2\nu_{\alpha}}} \right),
$$
then, for $k$ large enough, for all $\gamma\in \partial B_k$
$$
|\widetilde{f}(\gamma)-f_0(\gamma)|<|f_0(\gamma)|.
$$

Using Rouch\'{e}'s Theorem, we deduce that $\widetilde{f}$ and $f_0$ have the same number of zeros in $B_k$. Consequently, there exists a subsequence of roots of $\widetilde{f}$ that tends to the roots $\gamma_k^0$ of $f_0$, then there exists $N\in {\N}$ and a subsequence $\{\gamma_k\}_{|k|\geq N}$ of roots of $f(\gamma)$, such that $\gamma_k=\gamma_k^0+o(1)$ that tends to the roots $-\frac{2-\alpha}{2}\imath\left(k-\frac{\nu_{\alpha}}{2}+\frac{5}{4}\right)\pi$ of $f_0$. 
\\

Now, we can write
\begin{equation}\label{89}
\gamma_k=-\frac{2-\alpha}{2}\imath\left(k-\frac{\nu_{\alpha}}{2}+\frac{5}{4}\right)\pi+\varepsilon_k,
\end{equation}
then
\begin{eqnarray*}
e^{2\imath(z+(\nu_{\alpha}-1)\frac{\pi}{2}-\frac{\pi}{4})}&=&-e^{-\frac{4}{2-\alpha}\varepsilon_k}\\
&=&-1+\frac{4}{2-\alpha}\varepsilon_k+O(\varepsilon^2_k).
\end{eqnarray*}

Using the previous equation and the fact that $\widetilde{f}(\gamma_k)=0$, we get
\begin{eqnarray*}
\widetilde{f}(\gamma_k)&=&\frac{4}{2-\alpha}\varepsilon_k
-\frac{2\rho}{(1-\alpha)}\frac{c^-_{\nu_{\alpha,0}}}{c^+_{\nu_{\alpha,0}}}\left( \frac{2}{2-\alpha}\imath\right)^{-2\nu_{\alpha}}
\frac{\sin \nu_{\alpha}\pi}{{\left(-\frac{2-\alpha}{2}\imath k \pi \right)^{2\nu_{\alpha}}}}\\\\
&-& \Frac{(\frac{1}{2}-\nu_{\alpha})(\frac{3}{2}-\nu_{\alpha})}{2i}\left( \frac{2}{2-\alpha}\imath\right)^{-1}
\Frac{-2}{\left(-\frac{2-\alpha}{2}\imath k \pi \right)}\\ \\
&+& O(\varepsilon_k^2)+O\left(\frac{1}{k^{1+2\nu_{\alpha}}}\right)=0.\\ \\
\end{eqnarray*}
Hence
$$
\varepsilon_k=\Frac{\rho(2-\alpha)}{2(1-\alpha)}\frac{c^-_{\nu_{\alpha,0}}}{c^+_{\nu_{\alpha,0}}}\frac{\sin \nu_{\alpha}\pi}{(k\pi)^{2\nu_{\alpha}}}
+i\Frac{(2-\alpha)}{4}\Frac{(\frac{1}{2}-\nu_{\alpha})(\frac{3}{2}-\nu_{\alpha})}{k\pi}
+O\left(\frac{1}{k^{4\nu_{\alpha}}}\right),
$$
it follows that
$$
\gamma_k=-\frac{2-\alpha}{2}\imath\left(k-\frac{\nu_{\alpha}}{2}+\frac{5}{4}\right)\pi+
\Frac{\rho(2-\alpha)}{2(1-\alpha)}\frac{c^-_{\nu_{\alpha,0}}}{c^+_{\nu_{\alpha,0}}}\frac{\sin \nu_{\alpha}\pi}{(k\pi)^{2\nu_{\alpha}}}
+i\Frac{(2-\alpha)}{4}\Frac{(\frac{1}{2}-\nu_{\alpha})(\frac{3}{2}-\nu_{\alpha})}{k\pi}
+O\left(\frac{1}{k^{4\nu_{\alpha}}}\right).
$$
Since $\gamma_k^2=\imath\lambda_k$, then
\begin{eqnarray*}
\lambda_k&=&-\imath\gamma_k^2\\\\
&=&-\imath\left[- C_0^2(k\pi)^2-C_1^2\pi^2-2 C_0 C_1 k\pi^2-2{C_0 C_3}
+2i\Frac{C_0 C_2 \sin\nu_{\alpha}\pi}{(k\pi)^{2\nu_{\alpha}-1}}
+O\left(\frac{1}{k^{2\nu_{\alpha}}}\right)\right]\\\\
&=&i \left[C_0^2(k\pi)^2+2 C_0 C_1 k\pi^2\right]
+2\Frac{C_0 C_2 \sin\nu_{\alpha}\pi}{(k\pi)^{2\nu_{\alpha}-1}}+O(1)
\end{eqnarray*}
where
$$
C_0=-\frac{2-\alpha}{2}, C_1=-\frac{2-\alpha}{2}\left(-\frac{\nu_{\alpha}}{2}+\frac{5}{4}\right),
$$
$$
C_2=\Frac{\rho(2-\alpha)}{2(1-\alpha)}\frac{c^-_{\nu_{\alpha,0}}}{c^+_{\nu_{\alpha,0}}}
$$
and
$$
C_3=\Frac{(2-\alpha)}{4}{(\frac{1}{2}-\nu_{\alpha})(\frac{3}{2}-\nu_{\alpha})}.
$$

\hfill$\diamondsuit$\\

$\bullet\ \tilde{\alpha}>\Frac{4-3\alpha}{2(2-\alpha)}$.

We aim to solve the equation
$$
\matrix{f(\gamma)=-\imath\gamma(1-\alpha)\widetilde{d}^+
J_{1-\nu_\alpha}\left(\frac{2\gamma}{2-\alpha}\imath\right)+
\imath(1-\alpha)\rho(\lambda+\eta)^{\tilde{\alpha}-1}\widetilde{d}^-
J_{\nu_\alpha}\left(\frac{2\gamma}{2-\alpha}\imath\right)\hfill &\cr
+\rho \gamma(\lambda+\eta)^{\tilde{\alpha}-1}{\widetilde{d}}^-
J_{1+\nu_\alpha}\left(\frac{2\gamma}{2-\alpha}\imath\right)=0. &\cr}
$$
Using the development $(\ref{84})$, we get
$$
f(\gamma)=-\imath(1-\alpha)c^+_{\nu_{\alpha,0}}
\gamma^{1+\nu_{\alpha}}\left( \frac{2}{2-\alpha}\imath\right)^{\nu_{\alpha}}
\left( \frac{2}{\pi z}\right)^\frac{1}{2}
\frac{e^{-\imath(z+(\nu_{\alpha}-1)\frac{\pi}{2}-\frac{\pi}{4})}}{2}\widetilde{f}(\gamma),
$$
where $$z=\frac{2\gamma}{2-\alpha}\imath$$
and
\begin{eqnarray*}
\widetilde{f}(\gamma)&=& 1+e^{2\imath(z+(\nu_{\alpha}-1)\frac{\pi}{2}-\frac{\pi}{4})}-
\frac{\rho}{\imath}\frac{1}{1-\alpha}\frac{c^-_{\nu_{\alpha,0}}}{c^+_{\nu_{\alpha,0}}}\left( \frac{2}{2-\alpha}\imath\right)^{-2\nu_{\alpha}}(-i)^{\tilde{\alpha}-1}
\frac{e^{\imath(2z-3\frac{\pi}{2})}+e^{\imath\nu_{\alpha}\pi}}{\gamma^{2\nu_{\alpha}-2\tilde{\alpha}+2}}\\\\
&-&\Frac{(\frac{1}{2}-\nu_{\alpha})(\frac{3}{2}-\nu_{\alpha})}{2i}\left( \frac{2}{2-\alpha}\imath\right)^{-1}
\Frac{e^{2\imath(z+(\nu_{\alpha}-1)\frac{\pi}{2}-\frac{\pi}{4})}-1}{\gamma}\\ \\
&-&\rho\frac{1}{1-\alpha}\frac{c^-_{\nu_{\alpha,0}}}{c^+_{\nu_{\alpha,0}}}
\left( \frac{2}{2-\alpha}\imath\right) ^{-2\nu_{\alpha}-1}\Frac{(\nu_{\alpha}+\frac{1}{2})(\nu_{\alpha}+\frac{3}{2})}{2}
(-i)^{\tilde{\alpha}-1}\frac{e^{\imath(2z-3\frac{\pi}{2})}-e^{\imath\nu_{\alpha}\pi}}{\gamma^{3+2\nu_{\alpha}-2\tilde{\alpha}}}\\ \\
&-&\rho \frac{c^-_{\nu_{\alpha,0}}}{c^+_{\nu_{\alpha,0}}}
\left( \frac{2}{2-\alpha}\imath\right) ^{-2\nu_{\alpha}}\frac{e^{\imath(2z-\pi)}+e^{\imath(\nu_{\alpha}\pi-\frac{\pi}{2})}}{\gamma^{3+2\nu_{\alpha}-2\tilde{\alpha}}}+O\left(\frac{1}{\gamma^2} \right)\\\\
&=&f_0(\gamma)+\frac{f_1(\gamma)}{\gamma^{2\nu_{\alpha}-2\tilde{\alpha}+2}}++\frac{f_2(\gamma)}{\gamma}+\frac{f_3(\gamma)}{\gamma^{3+2\nu_{\alpha}-2\tilde{\alpha}}}+O\left(\frac{1}{\gamma^2} \right)
\end{eqnarray*}
with
$$
f_0(\gamma)=1+e^{2\imath(z+(\nu_{\alpha}-1)\frac{\pi}{2}-\frac{\pi}{4})},
$$
$$
f_1(\gamma)=-\frac{\rho}{\imath}\frac{1}{1-\alpha}\frac{c^-_{\nu_{\alpha,0}}}{c^+_{\nu_{\alpha,0}}}\left( \frac{2}{2-\alpha}\imath\right) ^{-2\nu_{\alpha}}
(e^{\imath(2z-3\frac{\pi}{2})}+e^{\imath\nu_{\alpha}\pi}).
$$
$$
f_2(\gamma)=-\Frac{(\frac{1}{2}-\nu_{\alpha})(\frac{3}{2}-\nu_{\alpha})}{2i}\left( \frac{2}{2-\alpha}\imath\right)^{-1}
(e^{2\imath(z+(\nu_{\alpha}-1)\frac{\pi}{2}-\frac{\pi}{4})}-1).
$$
and
$$
\matrix{f_3(\gamma)=&-&\rho\frac{1}{1-\alpha}\frac{c^-_{\nu_{\alpha,0}}}{c^+_{\nu_{\alpha,0}}}
\left( \frac{2}{2-\alpha}\imath\right) ^{-2\nu_{\alpha}-1}\Frac{(\nu_{\alpha}+\frac{1}{2})(\nu_{\alpha}+\frac{3}{2})}{2}
(-i)^{\tilde{\alpha}-1}(e^{\imath(2z-3\frac{\pi}{2})}-e^{\imath\nu_{\alpha}\pi}) \hfill \cr
&-&\rho \frac{c^-_{\nu_{\alpha,0}}}{c^+_{\nu_{\alpha,0}}}
\left( \frac{2}{2-\alpha}\imath\right) ^{-2\nu_{\alpha}}(e^{\imath(2z-\pi)}+e^{\imath(\nu_{\alpha}\pi-\frac{\pi}{2})})\hfill  \cr}
$$
Note that $f_0, f_1, f_2$ and $f_3$ remain bounded in the strip $-\alpha_0\leq\Re(\lambda)\leq0$.\\

We search the roots of $f_0$,
$$
f_0(\gamma)=0 \ \Leftrightarrow \ 1+e^{2\imath(z+(\nu_{\alpha}-1)\frac{\pi}{2}-\frac{\pi}{4})}=0,
$$
so, $f_0$ has the following roots
$$
\gamma_k^0=-\frac{2-\alpha}{2}\imath\left(k-\frac{\nu_{\alpha}}{2}+\frac{5}{4}\right)\pi, \ \ k \in \Z.
$$
Using Rouch\'{e}'s Theorem, we deduce that $\widetilde{f}$ and $f_0$ have the same number of zeros in $B_k$. Consequently, there exists a subsequence of roots of $\widetilde{f}$ that tends to the roots $\gamma_k^0$ of $f_0$, then there exists $N\in {\N}$ and a subsequence $\{\gamma_k\}_{|k|\geq N}$ of roots of $f(\gamma)$, such that $\gamma_k=\gamma_k^0+o(1)$ that tends to the roots $-\frac{2-\alpha}{2}\imath\left(k-\frac{\nu_{\alpha}}{2}+\frac{5}{4}\right)\pi$ of $f_0$. 
\\

Now, we can write
\begin{equation}\label{89n}
\gamma_k=-\frac{2-\alpha}{2}\imath\left(k-\frac{\nu_{\alpha}}{2}+\frac{5}{4}\right)\pi+\varepsilon_k,
\end{equation}
then
\begin{eqnarray*}
e^{2\imath(z+(\nu_{\alpha}-1)\frac{\pi}{2}-\frac{\pi}{4})}&=&-e^{-\frac{4}{2-\alpha}\varepsilon_k}\\
&=&-1+\frac{4}{2-\alpha}\varepsilon_k+O(\varepsilon^2_k).
\end{eqnarray*}
Using the previous equation and the fact that $\widetilde{f}(\gamma_k)=0$, we get
\begin{eqnarray*}
\widetilde{f}(\gamma_k)&=&\frac{4}{2-\alpha}\varepsilon_k
-\frac{2\rho}{(1-\alpha)}\frac{c^-_{\nu_{\alpha,0}}}{c^+_{\nu_{\alpha,0}}}\left( \frac{2}{2-\alpha}\imath\right)^{-2\nu_{\alpha}}
\frac{(-i)^{\tilde{\alpha}-1}\sin \nu_{\alpha}\pi}{{\left(-\frac{2-\alpha}{2}\imath k \pi \right)^{2\nu_{\alpha}-2\tilde{\alpha}+2}}}\\\\
&-& \Frac{(\frac{1}{2}-\nu_{\alpha})(\frac{3}{2}-\nu_{\alpha})}{2i}\left( \frac{2}{2-\alpha}\imath\right)^{-1}
\Frac{-2}{\left(-\frac{2-\alpha}{2}\imath k \pi \right)}
+ O(\varepsilon_k^2)+O\left(\frac{1}{k^{3+2\nu_{\alpha}-2\tilde{\alpha}}}\right)=0,\\ \\
\end{eqnarray*}
hence
$$
\varepsilon_k=\Frac{\rho(2-\alpha)}{2(1-\alpha)}\left(\Frac{2}{2-\alpha}\right)^{2-2\tilde{\alpha}}
\frac{c^-_{\nu_{\alpha,0}}}{c^+_{\nu_{\alpha,0}}}\frac{(-i)^{3\tilde{\alpha}-3}\sin \nu_{\alpha}\pi}{(k\pi)^{2\nu_{\alpha}-2\tilde{\alpha}+2}}
+i\Frac{(2-\alpha)}{4}\Frac{(\frac{1}{2}-\nu_{\alpha})(\frac{3}{2}-\nu_{\alpha})}{k\pi}
+O\left(\frac{1}{k^{2(2\nu_{\alpha}-2\tilde{\alpha}+2)}}\right),
$$
it follows that
$$
\matrix{\gamma_k&=&-\frac{2-\alpha}{2}\imath\left(k-\frac{\nu_{\alpha}}{2}+\frac{5}{4}\right)\pi+
\Frac{\rho(2-\alpha)}{2i(1-\alpha)}\left(\Frac{2}{2-\alpha}\right)^{2-2\tilde{\alpha}}
\frac{c^-_{\nu_{\alpha,0}}}{c^+_{\nu_{\alpha,0}}}\frac{(-i)^{3\tilde{\alpha}}\sin \nu_{\alpha}\pi}{(k\pi)^{2\nu_{\alpha}-2\tilde{\alpha}+2}}\hfill \cr
&&+i\Frac{(2-\alpha)}{4}\Frac{(\frac{1}{2}-\nu_{\alpha})(\frac{3}{2}-\nu_{\alpha})}{k\pi}
+O\left(\frac{1}{k^{2(2\nu_{\alpha}-2\tilde{\alpha}+2)}}\right).\hfill \cr}
$$
Since $\gamma_k^2=\imath\lambda_k$, then
\begin{eqnarray*}
\lambda_k&=&-\imath\gamma_k^2\\\\
&=&-\imath\left[- C_0^2(k\pi)^2-C_1^2\pi^2-2 C_0 C_1 k\pi^2-2{C_0 C_3}
+2\Frac{C_0 C_2(-i)^{3\tilde{\alpha}}\sin \nu_{\alpha}\pi}{(k\pi)^{2\nu_{\alpha}-2\tilde{\alpha}+1}}
+O\left(\frac{1}{k^{2\nu_{\alpha}-2\tilde{\alpha}+2}}\right)\right]\\\\
&=&i \left[C_0^2(k\pi)^2+2 C_0 C_1 k\pi^2\right]
-2i\Frac{C_0 C_2(-i)^{3\tilde{\alpha}}\sin\nu_{\alpha} \pi}{(k\pi)^{2\nu_{\alpha}-2\tilde{\alpha}+1}}+O\left(1\right),
\end{eqnarray*}
where
$$
C_0=-\frac{2-\alpha}{2}, C_1=-\frac{2-\alpha}{2}\left(-\frac{\nu_{\alpha}}{2}+\frac{5}{4}\right),
$$
$$
C_2=\Frac{\rho(2-\alpha)}{2(1-\alpha)}\left(\Frac{2}{2-\alpha}\right)^{2-2\tilde{\alpha}}\frac{c^-_{\nu_{\alpha,0}}}{c^+_{\nu_{\alpha,0}}}
$$
and
$$
C_3=\Frac{(2-\alpha)}{4}{(\frac{1}{2}-\nu_{\alpha})(\frac{3}{2}-\nu_{\alpha})}.
$$

\hfill$\diamondsuit$\\

$\bullet\ \nu_{\alpha}<\tilde{\alpha}<\Frac{4-3\alpha}{2(2-\alpha)}$.

We aim to solve the equation
$$
\matrix{f(\gamma)=-\imath\gamma(1-\alpha)\widetilde{d}^+
J_{1-\nu_\alpha}\left(\frac{2\gamma}{2-\alpha}\imath\right)+
\imath(1-\alpha)\rho(\lambda+\eta)^{\tilde{\alpha}-1}\widetilde{d}^-
J_{\nu_\alpha}\left(\frac{2\gamma}{2-\alpha}\imath\right)\hfill &\cr
+\rho \gamma(\lambda+\eta)^{\tilde{\alpha}-1}{\widetilde{d}}^-
J_{1+\nu_\alpha}\left(\frac{2\gamma}{2-\alpha}\imath\right)=0. &\cr}
$$
Using the development $(\ref{84})$, we get
$$
f(\gamma)=-\imath(1-\alpha)c^+_{\nu_{\alpha,0}}
\gamma^{1+\nu_{\alpha}}\left( \frac{2}{2-\alpha}\imath\right)^{\nu_{\alpha}}
\left( \frac{2}{\pi z}\right)^\frac{1}{2}
\frac{e^{-\imath(z+(\nu_{\alpha}-1)\frac{\pi}{2}-\frac{\pi}{4})}}{2}\widetilde{f}(\gamma),
$$
where $$z=\frac{2\gamma}{2-\alpha}\imath$$
and
\begin{eqnarray*}
\widetilde{f}(\gamma)&=& 1+e^{2\imath(z+(\nu_{\alpha}-1)\frac{\pi}{2}-\frac{\pi}{4})}
-\Frac{(\frac{1}{2}-\nu_{\alpha})(\frac{3}{2}-\nu_{\alpha})}{2i}\left( \frac{2}{2-\alpha}\imath\right)^{-1}
\Frac{e^{2\imath(z+(\nu_{\alpha}-1)\frac{\pi}{2}-\frac{\pi}{4})}-1}{\gamma}\\ \\
&-&\frac{\rho}{\imath}\frac{1}{1-\alpha}\frac{c^-_{\nu_{\alpha,0}}}{c^+_{\nu_{\alpha,0}}}\left( \frac{2}{2-\alpha}\imath\right)^{-2\nu_{\alpha}}(-i)^{\tilde{\alpha}-1}
\frac{e^{\imath(2z-3\frac{\pi}{2})}+e^{\imath\nu_{\alpha}\pi}}{\gamma^{2\nu_{\alpha}-2\tilde{\alpha}+2}}\\\\
&+& \Frac{(\frac{1}{2}-\nu_{\alpha})(\frac{3}{2}-\nu_{\alpha})(\frac{1}{2}+\nu_{\alpha})(\frac{5}{2}-\nu_{\alpha})}{8}
\left( \frac{2}{2-\alpha}\imath\right)^{-2}
\Frac{1+e^{2\imath(z+(\nu_{\alpha}-1)\frac{\pi}{2}-\frac{\pi}{4})}}{\gamma^2}\\\\
&-&\rho\frac{1}{1-\alpha}\frac{c^-_{\nu_{\alpha,0}}}{c^+_{\nu_{\alpha,0}}}
\left( \frac{2}{2-\alpha}\imath\right)^{-2\nu_{\alpha}-1}\Frac{(\nu_{\alpha}+\frac{1}{2})(\nu_{\alpha}+\frac{3}{2})}{2}
(-i)^{\tilde{\alpha}-1}\frac{e^{\imath(2z-3\frac{\pi}{2})}-e^{\imath\nu_{\alpha}\pi}}{\gamma^{3+2\nu_{\alpha}-2\tilde{\alpha}}}\\ \\
&-&\rho \frac{c^-_{\nu_{\alpha,0}}}{c^+_{\nu_{\alpha,0}}}
\left( \frac{2}{2-\alpha}\imath\right)^{-2\nu_{\alpha}}(-i)^{\tilde{\alpha}-1}\frac{e^{\imath(2z-\pi)}+e^{\imath(\nu_{\alpha}\pi-\frac{\pi}{2})}}{\gamma^{3+2\nu_{\alpha}-2\tilde{\alpha}}}+O\left(\frac{1}{\gamma^3} \right)\\\\
&=&f_0(\gamma)+\frac{f_1(\gamma)}{\gamma}+\frac{f_2(\gamma)}{\gamma^{2\nu_{\alpha}-2\tilde{\alpha}+2}}+
\frac{f_3(\gamma)}{\gamma^2}+\frac{f_4(\gamma)}{\gamma^{3+2\nu_{\alpha}-2\tilde{\alpha}}}+O\left(\frac{1}{\gamma^{3}} \right)
\end{eqnarray*}
with
$$
f_0(\gamma)=1+e^{2\imath(z+(\nu_{\alpha}-1)\frac{\pi}{2}-\frac{\pi}{4})},
$$
$$
f_1(\gamma)=-\Frac{(\frac{1}{2}-\nu_{\alpha})(\frac{3}{2}-\nu_{\alpha})}{2i}\left( \frac{2}{2-\alpha}\imath\right)^{-1}
(e^{2\imath(z+(\nu_{\alpha}-1)\frac{\pi}{2}-\frac{\pi}{4})}-1).
$$
$$
f_2(\gamma)=-\frac{\rho}{\imath}\frac{1}{1-\alpha}\frac{c^-_{\nu_{\alpha,0}}}{c^+_{\nu_{\alpha,0}}}\left( \frac{2}{2-\alpha}\imath\right) ^{-2\nu_{\alpha}}
(e^{\imath(2z-3\frac{\pi}{2})}+e^{\imath\nu_{\alpha}\pi}).
$$
$$
f_3(\gamma)=\Frac{(\frac{1}{2}-\nu_{\alpha})(\frac{3}{2}-\nu_{\alpha})(\frac{1}{2}+\nu_{\alpha})(\frac{5}{2}-\nu_{\alpha})}{8}
\left( \frac{2}{2-\alpha}\imath\right)^{-2}
\Frac{1+e^{2\imath(z+(\nu_{\alpha}-1)\frac{\pi}{2}-\frac{\pi}{4})}}{\gamma^2}
$$
and
$$
\matrix{f_4(\gamma)=
&-&\rho\frac{1}{1-\alpha}\frac{c^-_{\nu_{\alpha,0}}}{c^+_{\nu_{\alpha,0}}}
\left( \frac{2}{2-\alpha}\imath\right)^{-2\nu_{\alpha}-1}\Frac{(\nu_{\alpha}+\frac{1}{2})(\nu_{\alpha}+\frac{3}{2})}{2}
(-i)^{\tilde{\alpha}-1}(e^{\imath(2z-3\frac{\pi}{2})}-e^{\imath\nu_{\alpha}\pi})\hfill \cr
&&-\rho \frac{c^-_{\nu_{\alpha,0}}}{c^+_{\nu_{\alpha,0}}}
\left( \frac{2}{2-\alpha}\imath\right)^{-2\nu_{\alpha}}(-i)^{\tilde{\alpha}-1}(e^{\imath(2z-\pi)}+e^{\imath(\nu_{\alpha}\pi-\frac{\pi}{2})})\hfill \cr}
$$
Note that $f_0$, $f_1$ and $f_2$ remain bounded in the strip $-\alpha_0\leq\Re(\lambda)\leq0$.\\

We search the roots of $f_0$,
$$
f_0(\gamma)=0 \ \Leftrightarrow \ 1+e^{2\imath(z+(\nu_{\alpha}-1)\frac{\pi}{2}-\frac{\pi}{4})}=0,
$$
so, $f_0$ has the following roots
$$
\gamma_k^0=-\frac{2-\alpha}{2}\imath\left(k-\frac{\nu_{\alpha}}{2}+\frac{5}{4}\right)\pi, \ \ k \in \Z.
$$
Using Rouch\'{e}'s Theorem, we deduce that $\widetilde{f}$ and $f_0$ have the same number of zeros in $B_k$. Consequently, there exists a subsequence of roots of $\widetilde{f}$ that tends to the roots $\gamma_k^0$ of $f_0$, then there exists $N\in {\N}$ and a subsequence $\{\gamma_k\}_{|k|\geq N}$ of roots of $f(\gamma)$, such that $\gamma_k=\gamma_k^0+o(1)$ that tends to the roots $-\frac{2-\alpha}{2}\imath\left(k-\frac{\nu_{\alpha}}{2}+\frac{5}{4}\right)\pi$ of $f_0$. 
\\

Now, we can write
\begin{equation}\label{89nn}
\gamma_k=-\frac{2-\alpha}{2}\imath\left(k-\frac{\nu_{\alpha}}{2}+\frac{5}{4}\right)\pi+\varepsilon_k,
\end{equation}
then
\begin{eqnarray*}
e^{2\imath(z+(\nu_{\alpha}-1)\frac{\pi}{2}-\frac{\pi}{4})}&=&-e^{-\frac{4}{2-\alpha}\varepsilon_k}\\
&=&-1+\frac{4}{2-\alpha}\varepsilon_k+O(\varepsilon^2_k).
\end{eqnarray*}
Using the previous equation and the fact that $\widetilde{f}(\gamma_k)=0$, we get
\begin{eqnarray*}
\widetilde{f}(\gamma_k)&=&\frac{4}{2-\alpha}\varepsilon_k
-\frac{2\rho}{(1-\alpha)}\frac{c^-_{\nu_{\alpha,0}}}{c^+_{\nu_{\alpha,0}}}\left( \frac{2}{2-\alpha}\imath\right)^{-2\nu_{\alpha}}\frac{(-i)^{\tilde{\alpha}-1}\sin\nu_{\alpha} \pi}{{\left(-\frac{2-\alpha}{2}\imath k \pi \right)^{2\nu_{\alpha}-2\tilde{\alpha}+2}}}\\\\
&-& \Frac{(\frac{1}{2}-\nu_{\alpha})(\frac{3}{2}-\nu_{\alpha})}{2i}\left( \frac{2}{2-\alpha}\imath\right)^{-1}
\Frac{-2}{\left(-\frac{2-\alpha}{2}\imath k \pi \right)}
+ O(\varepsilon_k^2)+O\left(\frac{1}{k^{2\nu_{\alpha}-2\tilde{\alpha}+3}}\right)=0.\\ \\
\end{eqnarray*}
Hence
$$
\varepsilon_k=i\Frac{(2-\alpha)}{4}\Frac{(\frac{1}{2}-\nu_{\alpha})(\frac{3}{2}-\nu_{\alpha})}{k\pi}+
\Frac{\rho(2-\alpha)}{2(1-\alpha)}\left(\Frac{2}{2-\alpha}\right)^{2-2\tilde{\alpha}}
\frac{c^-_{\nu_{\alpha,0}}}{c^+_{\nu_{\alpha,0}}}\frac{(-i)^{3\tilde{\alpha}-3}\sin\nu_{\alpha} \pi}{(k\pi)^{2\nu_{\alpha}-2\tilde{\alpha}+2}}
+O\left(\frac{1}{k^{2\nu_{\alpha}-2\tilde{\alpha}+3}}\right),
$$
it follows that
$$
\matrix{\gamma_k&=&-\frac{2-\alpha}{2}\imath\left(k-\frac{\nu_{\alpha}}{2}+\frac{5}{4}\right)\pi
+i\Frac{(2-\alpha)}{4}\Frac{(\frac{1}{2}-\nu_{\alpha})(\frac{3}{2}-\nu_{\alpha})}{k\pi}\hfill \cr
&&+\Frac{\rho(2-\alpha)}{2i(1-\alpha)}\left(\Frac{2}{2-\alpha}\right)^{2-2\tilde{\alpha}}
\frac{c^-_{\nu_{\alpha,0}}}{c^+_{\nu_{\alpha,0}}}\frac{(-i)^{3\tilde{\alpha}}\sin\nu_{\alpha} \pi}{(k\pi)^{2\nu_{\alpha}-2\tilde{\alpha}+2}}
+O\left(\frac{1}{k^{2\nu_{\alpha}-2\tilde{\alpha}+3}}\right).\hfill \cr}
$$
Since $\gamma_k^2=\imath\lambda_k$, then
\begin{eqnarray*}
\lambda_k&=&-\imath\gamma_k^2\\\\
&=&-\imath\left[- C_0^2(k\pi)^2-C_1^2\pi^2-2 C_0 C_1 k\pi^2-2{C_0 C_3}
+2\Frac{C_0 C_2(-i)^{3\tilde{\alpha}}\sin\nu_{\alpha} \pi}{(k\pi)^{2\nu_{\alpha}-2\tilde{\alpha}+1}}
+O\left(\frac{1}{k}\right)\right]\\\\
&=&i \left[C_0^2(k\pi)^2+C_1^2\pi^2+2 C_0 C_1 k\pi^2+2{C_0 C_3}\right]
-2i\Frac{C_0 C_2(-i)^{3\tilde{\alpha}}\sin\nu_{\alpha} \pi}{(k\pi)^{2\nu_{\alpha}-2\tilde{\alpha}+1}}+O\left(\frac{1}{k}\right)
\end{eqnarray*}
where
$$
C_0=-\frac{2-\alpha}{2}, C_1=-\frac{2-\alpha}{2}\left(-\frac{\nu_{\alpha}}{2}+\frac{5}{4}\right),
$$
$$
C_2=\Frac{\rho(2-\alpha)}{2(1-\alpha)}\left(\Frac{2}{2-\alpha}\right)^{2-2\tilde{\alpha}}\frac{c^-_{\nu_{\alpha,0}}}{c^+_{\nu_{\alpha,0}}}
$$
and
$$
C_3=\Frac{(2-\alpha)}{4}{(\frac{1}{2}-\nu_{\alpha})(\frac{3}{2}-\nu_{\alpha})}.
$$

\hfill$\diamondsuit$\\

$\bullet\  \tilde{\alpha}<\nu_{\alpha}$.

We aim to solve the equation
$$
\matrix{f(\gamma)=-\imath\gamma(1-\alpha)\widetilde{d}^+
J_{1-\nu_\alpha}\left(\frac{2\gamma}{2-\alpha}\imath\right)+
\imath(1-\alpha)\rho(\lambda+\eta)^{\tilde{\alpha}-1}\widetilde{d}^-
J_{\nu_\alpha}\left(\frac{2\gamma}{2-\alpha}\imath\right)\hfill &\cr
+\rho \gamma(\lambda+\eta)^{\tilde{\alpha}-1}{\widetilde{d}}^-
J_{1+\nu_\alpha}\left(\frac{2\gamma}{2-\alpha}\imath\right)=0. &\cr}
$$
Using the development $(\ref{84})$, we get
$$
f(\gamma)=-\imath(1-\alpha)c^+_{\nu_{\alpha,0}}
\gamma^{1+\nu_{\alpha}}\left( \frac{2}{2-\alpha}\imath\right)^{\nu_{\alpha}}
\left( \frac{2}{\pi z}\right)^\frac{1}{2}
\frac{e^{-\imath(z+(\nu_{\alpha}-1)\frac{\pi}{2}-\frac{\pi}{4})}}{2}\widetilde{f}(\gamma),
$$
where $$z=\frac{2\gamma}{2-\alpha}\imath$$
and
\begin{eqnarray*}
\widetilde{f}(\gamma)&=& 1+e^{2\imath(z+(\nu_{\alpha}-1)\frac{\pi}{2}-\frac{\pi}{4})}
-\Frac{(\frac{1}{2}-\nu_{\alpha})(\frac{3}{2}-\nu_{\alpha})}{2i}\left( \frac{2}{2-\alpha}\imath\right)^{-1}
\Frac{e^{2\imath(z+(\nu_{\alpha}-1)\frac{\pi}{2}-\frac{\pi}{4})}-1}{\gamma}\\ \\
&+&\Frac{(\frac{1}{2}-\nu_{\alpha})(\frac{3}{2}-\nu_{\alpha})(\frac{1}{2}+\nu_{\alpha})(\frac{5}{2}-\nu_{\alpha})}{8}
\left( \frac{2}{2-\alpha}\imath\right)^{-2}
\Frac{1+e^{2\imath(z+(\nu_{\alpha}-1)\frac{\pi}{2}-\frac{\pi}{4})}}{\gamma^2}\\\\
&-&\frac{\rho}{\imath}\frac{1}{1-\alpha}\frac{c^-_{\nu_{\alpha,0}}}{c^+_{\nu_{\alpha,0}}}\left( \frac{2}{2-\alpha}\imath\right)^{-2\nu_{\alpha}}(-i)^{\tilde{\alpha}-1}
\frac{e^{\imath(2z-3\frac{\pi}{2})}+e^{\imath\nu_{\alpha}\pi}}{\gamma^{2\nu_{\alpha}-2\tilde{\alpha}+2}}
+O\left(\frac{1}{\gamma^3} \right)\\ \\
&=&f_0(\gamma)+\frac{f_1(\gamma)}{\gamma}+\frac{f_2(\gamma)}{\gamma^2}+\frac{f_3(\gamma)}{\gamma^{2\nu_{\alpha}-2\tilde{\alpha}+2}}
+O\left(\frac{1}{\gamma^3} \right)
\end{eqnarray*}
with
$$
f_0(\gamma)=1+e^{2\imath(z+(\nu_{\alpha}-1)\frac{\pi}{2}-\frac{\pi}{4})},
$$
$$
f_1(\gamma)=-\Frac{(\frac{1}{2}-\nu_{\alpha})(\frac{3}{2}-\nu_{\alpha})}{2i}\left( \frac{2}{2-\alpha}\imath\right)^{-1}
(e^{2\imath(z+(\nu_{\alpha}-1)\frac{\pi}{2}-\frac{\pi}{4})}-1)
$$

$$
f_2(\gamma)=\Frac{(\frac{1}{2}-\nu_{\alpha})(\frac{3}{2}-\nu_{\alpha})(\frac{1}{2}+\nu_{\alpha})(\frac{5}{2}-\nu_{\alpha})}{8}
\left( \frac{2}{2-\alpha}\imath\right)^{-2}
(1+e^{2\imath(z+(\nu_{\alpha}-1)\frac{\pi}{2}-\frac{\pi}{4})}).
$$
and
$$
f_3(\gamma)=-\frac{\rho}{\imath}\frac{1}{1-\alpha}\frac{c^-_{\nu_{\alpha,0}}}{c^+_{\nu_{\alpha,0}}}\left( \frac{2}{2-\alpha}\imath\right) ^{-2\nu_{\alpha}}
(-i)^{\tilde{\alpha}-1}(e^{\imath(2z-3\frac{\pi}{2})}+e^{\imath\nu_{\alpha}\pi}).
$$
Note that $f_0$, $f_1$ and $f_2$ remain bounded in the strip $-\alpha_0\leq\Re(\lambda)\leq0$.\\

We search the roots of $f_0$,
$$
f_0(\gamma)=0 \ \Leftrightarrow \ 1+e^{2\imath(z+(\nu_{\alpha}-1)\frac{\pi}{2}-\frac{\pi}{4})}=0,
$$
so, $f_0$ has the following roots
$$
\gamma_k^0=-\frac{2-\alpha}{2}\imath\left(k-\frac{\nu_{\alpha}}{2}+\frac{5}{4}\right)\pi, \ \ k \in \Z.
$$
Using Rouch\'{e}'s Theorem, we deduce that $\widetilde{f}$ and $f_0$ have the same number of zeros in $B_k$. Consequently, there exists a subsequence of roots of $\widetilde{f}$ that tends to the roots $\gamma_k^0$ of $f_0$, then there exists $N\in {\N}$ and a subsequence $\{\gamma_k\}_{|k|\geq N}$ of roots of $f(\gamma)$, such that $\gamma_k=\gamma_k^0+o(1)$ that tends to the roots $-\frac{2-\alpha}{2}\imath\left(k-\frac{\nu_{\alpha}}{2}+\frac{5}{4}\right)\pi$ of $f_0$. 
\\

Now, we can write
\begin{equation}\label{89n4}
\gamma_k=-\frac{2-\alpha}{2}\imath\left(k-\frac{\nu_{\alpha}}{2}+\frac{5}{4}\right)\pi+\varepsilon_k,
\end{equation}
then
\begin{eqnarray*}
e^{2\imath(z+(\nu_{\alpha}-1)\frac{\pi}{2}-\frac{\pi}{4})}&=&-e^{-\frac{4}{2-\alpha}\varepsilon_k}\\
&=&-1+\frac{4}{2-\alpha}\varepsilon_k+O(\varepsilon^2_k).
\end{eqnarray*}

Using the previous equation and the fact that $\widetilde{f}(\gamma_k)=0$, we get
$$
\widetilde{f}(\gamma_k)=\frac{4}{2-\alpha}\varepsilon_k
- \Frac{(\frac{1}{2}-\nu_{\alpha})(\frac{3}{2}-\nu_{\alpha})}{2i}\left( \frac{2}{2-\alpha}\imath\right)^{-1}
\Frac{-2}{\left(-\frac{2-\alpha}{2}\imath k \pi \right)}+o(\varepsilon_k)+o\left(\Frac{1}{k}\right)=0.
$$
hence
$$
\varepsilon_k=
i\Frac{(2-\alpha)}{4}\Frac{(\frac{1}{2}-\nu_{\alpha})(\frac{3}{2}-\nu_{\alpha})}{k\pi}
+o(\varepsilon_k)+o\left(\Frac{1}{k}\right).
$$
We can write
\begin{equation}\label{332}
\gamma_k=-\frac{2-\alpha}{2}\imath\left(k-\frac{\nu_{\alpha}}{2}+\frac{5}{4}\right)\pi+i\Frac{(2-\alpha)}{4}\Frac{(\frac{1}{2}-\nu_{\alpha})(\frac{3}{2}-\nu_{\alpha})}{k\pi}
+\tilde{\varepsilon}_k,
\end{equation}
where $\tilde{\varepsilon}_k= o\left(\Frac{1}{k}\right)$.
$$
\matrix{e^{2\imath(z+(\nu_{\alpha}-1)\frac{\pi}{2}-\frac{\pi}{4})}&=&-e^{-\frac{4}{2-\alpha}(\frac{il}{k\pi}+\tilde{\varepsilon}_k)}\hfill \cr
&=&-1+\frac{4}{2-\alpha}(\frac{il}{k\pi}+\tilde{\varepsilon}_k)
-\Frac{1}{2}\frac{4}{2-\alpha}\left(\frac{il}{k\pi}\right)^2
+O(\tilde{\varepsilon}_k^2)+o\left(\Frac{1}{k^2}\right)+O\left(\Frac{\tilde{\varepsilon}_k}{k}\right).\hfill \cr}
$$
where $l=\Frac{(2-\alpha){(1/2-\nu_{\alpha})(3/2-\nu_{\alpha})}}{4}$.
Substituting (\ref{332}) into $\widetilde{f}(\gamma_k)=0$,
we get
$$
\widetilde{f}(\gamma_k)=\frac{4}{2-\alpha}\varepsilon_k-i\Frac{2m (-\frac{\nu_{\alpha}}{2}+\frac{5}{4})}{k^2\pi}+
o\left(\Frac{1}{k^2}\right)+O(\tilde{\varepsilon}_k^2).
$$
where $m=-\Frac{{(1/2-\nu_{\alpha})(3/2-\nu_{\alpha})}}{2}$.
hence
$$
\tilde{\varepsilon}_k=-i\frac{2-\alpha}{2}\Frac{m (-\frac{\nu_{\alpha}}{2}+\frac{5}{4})}{k^2\pi}
+o\left(\Frac{1}{k^2}\right)+O(\tilde{\varepsilon}_k^2).
$$
We can write
\begin{equation}\label{333}
\gamma_k=-\frac{2-\alpha}{2}\imath\left(k-\frac{\nu_{\alpha}}{2}+\frac{5}{4}\right)\pi+i\Frac{(2-\alpha)}{4}\Frac{(\frac{1}{2}-\nu_{\alpha})(\frac{3}{2}-\nu_{\alpha})}{k\pi}
-i\frac{2-\alpha}{2}\Frac{m (-\frac{\nu_{\alpha}}{2}+\frac{5}{4})}{k^2\pi}+\tilde{\tilde\varepsilon}_k,
\end{equation}
where $\tilde{\tilde\varepsilon}_k= o\left(\Frac{1}{k^2}\right)$.
Substituting (\ref{333}) into $\widetilde{f}(\gamma_k)=0$,
we get
$$
\widetilde{f}(\gamma_k)=\frac{4}{2-\alpha}\tilde{\tilde\varepsilon}_k
-\frac{2\rho}{(1-\alpha)}\frac{c^-_{\nu_{\alpha,0}}}{c^+_{\nu_{\alpha,0}}}\left( \frac{2}{2-\alpha}\imath\right)^{-2\nu_{\alpha}}\frac{(-i)^{\tilde{\alpha}-1}\sin\nu_{\alpha} \pi}{{\left(-\frac{2-\alpha}{2}\imath k \pi \right)^{2\nu_{\alpha}-2\tilde{\alpha}+2}}}
+ O(\tilde{\tilde\varepsilon}_k^2)+O\left(\frac{1}{k^{3}}\right).
$$
hence
$$
\tilde{\tilde\varepsilon}_k=\Frac{\rho(2-\alpha)}{2(1-\alpha)}\left(\Frac{2}{2-\alpha}\right)^{2-2\tilde{\alpha}}\frac{c^-_{\nu_{\alpha,0}}}{c^+_{\nu_{\alpha,0}}}
\frac{(-i)^{3\tilde{\alpha}-3}\sin\nu_{\alpha} \pi}{(k\pi)^{2\nu_{\alpha}-2\tilde{\alpha}+2}}
+O\left(\frac{1}{k^{3}}\right).
$$
it follows that
$$
\matrix{\gamma_k=
-\frac{2-\alpha}{2}\imath\left(k-\frac{\nu_{\alpha}}{2}+\frac{5}{4}\right)\pi
+i\Frac{(2-\alpha)}{4}\Frac{(\frac{1}{2}-\nu_{\alpha})(\frac{3}{2}-\nu_{\alpha})}{k\pi}\hfill &\cr
\hspace{1cm}-i\frac{2-\alpha}{2}\Frac{m (-\frac{\nu_{\alpha}}{2}+\frac{5}{4})}{k^2\pi}
+\Frac{\rho(2-\alpha)}{2i(1-\alpha)}\left(\Frac{2}{2-\alpha}\right)^{2-2\tilde{\alpha}}
\frac{c^-_{\nu_{\alpha,0}}}{c^+_{\nu_{\alpha,0}}}\frac{(-i)^{3\tilde{\alpha}}\sin\nu_{\alpha} \pi}{(k\pi)^{2\nu_{\alpha}-2\tilde{\alpha}+2}}
+O\left(\frac{1}{k^3}\right)\hfill &\cr}
$$
Since $\gamma_k^2=\imath\lambda_k$, then
$$
\matrix{\lambda_k&=&-\imath\gamma_k^2 \hfill \cr
&=&-\imath\left[- C_0^2(k\pi)^2-C_1^2\pi^2-2 C_0 C_1 k\pi^2-2{C_0 C_2}-2\frac{C_0 C_3}{k}-2\frac{C_0 C_1}{k}
+2\Frac{C_0 C_4(-i)^{3\tilde{\alpha}}\sin\nu_{\alpha} \pi}{(k\pi)^{2\nu_{\alpha}-2\tilde{\alpha}+1}}\right.\hfill \cr
&&\left.+O\left(\frac{1}{k^{2}}\right)\right] \cr
&=&i \left[C_0^2(k\pi)^2+C_1^2\pi^2+2 C_0 C_1 k\pi^2+2{C_0 C_2}+2\frac{C_0 C_3}{k}+2\frac{C_0 C_1}{k}\right]
-2i\Frac{C_0 C_4(-i)^{3\tilde{\alpha}}\sin\nu_{\alpha} \pi}{(k\pi)^{2\nu_{\alpha}-2\tilde{\alpha}+1}}\hfill \cr
&&+O\left(\frac{1}{k^{2}}\right)  \cr}
$$
where
$$
C_0=-\frac{2-\alpha}{2}, C_1=-\frac{2-\alpha}{2}\left(-\frac{\nu_{\alpha}}{2}+\frac{5}{4}\right),
$$
$$
C_2=\Frac{(2-\alpha)}{4}(\frac{1}{2}-\nu_{\alpha})(\frac{3}{2}-\nu_{\alpha}),
$$
and
$$
C_3=-\frac{2-\alpha}{2}m (-\frac{\nu_{\alpha}}{2}+\frac{5}{4})
$$
$$
C_4=\Frac{\rho(2-\alpha)}{2(1-\alpha)}\left(\Frac{2}{2-\alpha}\right)^{2-2\tilde{\alpha}}\frac{c^-_{\nu_{\alpha,0}}}{c^+_{\nu_{\alpha,0}}}.
$$
Now, setting ${\tilde V}_k=(\lambda_k^0 I-{\cal A}) V_k$, where $V_k$ is
a normalized eigenfunction associated to $\lambda_k$ and
$$
\lambda_k^0=\left\{\matrix{i \left[C_0^2(k\pi)^2+C_1^2\pi^2+2 C_0 C_1 k\pi^2+2{C_0 C_3}\right]\hfill & \hbox{ if }
\nu_{\alpha} <\tilde{\alpha}<\Frac{4-3\alpha}{2(2-\alpha)},\hfill \cr
i \left[C_0^2(k\pi)^2+C_1^2\pi^2+2 C_0 C_1 k\pi^2+2{C_0 C_2}+2\frac{C_0 C_3}{k}+2\frac{C_0 C_1}{k}\right]\hfill &
\hbox{ if } \tilde{\alpha}<\nu_{\alpha}.\hfill \cr}\right.
$$
We then have
$$
\matrix{\|(\lambda_k^0 I-{\cal A})^{-1}\|_{{\cal L}({\cal H})}=\Sup_{V\in {\cal H}, V\not=0} \Frac{\|(\lambda_k^0 I-{\cal A})^{-1}V\|_{{\cal H}}}{\|V\|_{{\cal H}}}
&\geq &\Frac{\|(\lambda_k^0 I-{\cal A})^{-1}{\tilde V}_k\|_{{\cal H}}}{\|{\tilde V}_k\|_{{\cal H}}}\hfill \cr
&\geq &\Frac{\|V_k\|_{{\cal H}}}{\|(\lambda_k^0 I-{\cal A}) V_k\|_{{\cal H}}}.\hfill \cr}
$$
Hence, by Lemma \ref{l4.1}, we deduce that
$$
\|(\lambda_k^0 I-{\cal A})^{-1}\|_{{\cal L}({\cal H})}\geq
c |k|^{2\nu_{\alpha}-2\tilde{\alpha}+1}\equiv |\lambda_k^0|^{\nu_{\alpha}-\tilde{\alpha}+1/2}\quad\hbox{ if } 0< \tilde{\alpha}<\Frac{4-3\alpha}{2(2-\alpha)},
$$
Thus, the second condition of Lemma \ref{4.1} is not satisfied for $0< \tilde{\alpha}<\Frac{4-3\alpha}{2(2-\alpha)}$. So that, the semigroup
$e^{t {\cal A}}$ is not exponentially stable. Thus the proof is complete.

\hfill$\diamondsuit$\\

\subsection*{Conclusion}
The problem (\ref{P}) exhibits strong degeneracy at zero, which leads to exponential
or polynomial stabilization of the system depending on the relation between $\alpha$ and $\tilde{\alpha}$.
On the other hand, when the weight is non-degenerate or the damping is acting on the non-degenerate point $x=1$ one
can expect decay estimate similar to the case $\alpha=0$.
In this case, multiplier methods can be effectively used to derive energy estimates and demonstrate decay of solutions.

Here, we obtain a strong asymptotic behavior for $\eta\geq 0$ and when $\eta> 0$, we get sharp estimate for the
rate of energy decay of classical solutions depending on parameters $\alpha$ and $\tilde{\alpha}$. Our approach
is based on the asymptotic theory of $C_0$- semigroups and in particular on a result due to Borichev and Tomilov
{\bf\cite{BT}}, which reduces the problem of estimating the rate of energy decay
to finding a growth bound for the resolvent of the semigroup generator by using Bessel functions.
In particular, we obtain uniform decay estimates for a weakly Schr\"{o}dinger equation under a weak damping.

\section*{Appendix A. Proof of Lemma \ref{l5.3}}

\noindent
We will use the following result.
\begin{lemma}\label{6.1}
For $a$ a complex number and $\Re{\nu_{\alpha}}>-1$, we have
$$
2a^2\int_{0}^{x}t(J_{\nu_{\alpha}}(at))^2dt=(a^2x^2-{\nu_{\alpha}}^2)(J_{\nu_{\alpha}}(ax))^2+\left(x\frac{d}{dx}(J_{\nu_{\alpha}}(ax))\right)^2.
$$
\end{lemma}
\noindent
{\bf Proof.}
We have
\begin{equation}\label{113}
\|\theta_+\|_{L^2(0,1)}^2=\int_{0}^{1}x^{1-\alpha}\left(J_{\nu_\alpha}\left( \frac{2}{2-\alpha}\mu x^{\frac{2-\alpha}{2}}\right)  \right)^2 dx.
\end{equation}

Suppose that $s=\frac{2}{2-\alpha}\mu x^{\frac{2-\alpha}{2}}$ in $(\ref{113})$, we find
$$
\|\theta_+\|_{L^2(0,1)}^2=\frac{2-\alpha}{2\mu^2}\int_{0}^{r}s(J_{\nu_\alpha}(s))^2ds, \hbox{ with } r=\frac{2\mu}{2-\alpha},
$$
using lemma $(\ref{6.1})$, we get
$$
\|\theta_+\|_{L^2(0,1)}^2=\frac{1}{2-\alpha}\frac{1}{r^2}
\left[ (r^2-\nu_\alpha^2)(J_{\nu_\alpha}(r))^2+(rj'_{\nu_\alpha}(r))^2\right],
$$
the relation $(\ref{45})$ gives
$$
\|\theta_+\|_{L^2(0,1)}^2=\frac{1}{2-\alpha}\frac{1}{r^2}
\left[(rJ_{\nu_\alpha}(r))^2+(rJ_{\nu_\alpha+1}(r))^2-2\nu_\alpha rJ_{\nu_\alpha}(r)J_{\nu_\alpha+1}(r)\right].
$$
In a similay way we prove the other inequalities.

\subsection*{Acknowledgments}
The authors would like to thank very much the referee for their constructive comments and suggestions
that helped to improve this article.


\begin{thebibliography}{99}
\bibitem{achour.1} Achouri, Z., Amroun, N., Benaissa, A.: {The Euler-Bernoulli beam equation with boundary
dissipation of fractional derivative type\/,} {Mathematical Methods in the Applied Sciences.} {\bf 40} (2017)-11, 3887-3854.
\bibitem{AB} Arendt, W., Batty, C.J.K.:{Tauberian theorems and stability of one-parameter semigroups,}
Trans Amer Math Soc.  {\bf 306} (1988)-2, 837-852.
\bibitem{BT} Borichev, A., Tomilov, Y.: {Optimal polynomial decay of functions and operator semigroups,}
Math Ann. {\bf 347} (2010)-2, 455-478.
\bibitem{Bh} Br\'{e}zis, H.: {Operateurs Maximaux Monotones et semi-groupes de contractions dans les espaces de Hilbert,}
Notas de matem\`{a}tica (50). North-Holland, Amsterdam: Universidade Federal do Rio de Janeiro and
University of Rochester; 1973.
\bibitem{choi} Choi, J. U., Maccamy, R. C.: {Fractional order Volterra equations with applications to elasticity\/,}
{J.Math. Anal. Appl.} {\bf 139} (1989), 448-464.
\bibitem{D} Dieudonn\'e, J.: {Calcul infinit\'{e}simal,} Collection Methodes, Herman, Paris. (1968).
\bibitem{jawad}  Fragnelli, G. Moumni, A. Salhi, J.: {Controllability and stabilization of a degenerate/singular
 Schr\"{o}dinger equation,} J. Math. Anal. Appl. {\bf 537} (2024)-2, Paper No. 128290, 32 pp.
\bibitem{KZ} Komornik, V.,Zuazua, E.: {A direct method for the boundary stabilization of the wave equation,}
J. Math. pures et appl., {\bf 69} (1990),33-54.
\bibitem{KGS} Krstic, M., Guo, B. Z., Smyshlyaev, A.:, {Boundary Controllers and Observers for the Linearized
Schr\"{o}dinger Equation,} SIAM Journal on Control and Optimization,  {\bf 49}, (2011)-4, 1479-1497.
\bibitem{Ire-Trig} Lasiecka, I., Triggiani, R.: {Well-posedness and sharp uniform decay rates at the
$L^2(\Omega)$-Level of the Schr\"{o}dinger equation with nonlinear boundary dissipation\/,}
J. evol. equ. {\bf 6} (2006), 485-537.
\bibitem{lebeau} Lebeau, G.: {Contr\^ole de l'\'equation de Schr\"{o}dinger\/,} J. Math. Pures Appl. {\bf 71} (1992) 267-291.
\bibitem{Le} Lebedev, N. N.: {Special Functions and their Applications,} New York: Dover Publications; 1972.
\bibitem{LGM} Luo,  Z. H., Guo  B. Z., Morgul, O.: {Stability  and  Stabilization  of  Infinite  Dimensional
Systems  with Applications,} Communications and Control Engineering Series. London: Springer; 1999.
\bibitem{M}  Machtyngier, E.: {Exact controllability for the Sch\"{o}rdinger equation,} SIAM J. Control Optim.
{\bf 32} (1994)-1, 24-34.
\bibitem{My-Zau} Machtyngier, E., Zuazua, E.: {Stabilization of the Schr\"{o}dinger equation\/,}
Portugaliae Mathematica. {\bf 51} (1994)-2, 243-256.
\bibitem{mbod} Mbodje, B.: {Wave energy decay under fractional derivative controls\/,} {IMA Journal of Mathematical
Control and Information.,}  {\bf 23} (2006), 237-257.
\bibitem{P} Pr\"{u}ss, J.:  {On  the  spectrum  of  C$_0$-semigroups,} Transactions  of  the  American
Mathematical  Society. {\bf 284} (1984)-2, 847-857.
\bibitem{Rebiai} Nicaise, S.,Rebiai, S.: {\em Stabilization of the Schr\"{o}dinger equation with a
delay term in boundary feedback or internal feedback\/,} Portugaliae Mathematica, {\bf 68} (2011)-1, 19-39.
\bibitem{W} Watson,  G. N.: {A  treatise  on  the  theory  of  Bessel  functions,  second  edition,}
Cambridge,  England: Cambridge University Press; 1944.
\bibitem{ZBA} Zerkouk, H., Aichi, C., Benaissa, A.:, {On the Stability of a Degenerate Wave Equation Under
Fractional Feedbacks Acting on the Degenerate Boundary,} J Dyn Control Syst, {\bf 28} (2022), 601-633.










\end{thebibliography}
\end{document}